\documentclass[11pt]{article}
\usepackage[margin=1in]{geometry}
\usepackage[T1]{fontenc}
\usepackage{lmodern}

\usepackage{amsmath}
\usepackage{amssymb}
\usepackage{amsthm}
\usepackage[numbers,sort&compress]{natbib}
\usepackage[title]{appendix}
\usepackage{mathtools}
\usepackage{graphicx}
\usepackage{booktabs}
\usepackage{array}
\newcolumntype{L}[1]{>{\raggedright\arraybackslash}p{#1}}
\usepackage{tikz}
\usepackage{pgfplots}
\pgfplotsset{compat=1.18}
\usepackage{adjustbox}
\usepackage{microtype}
\usepackage{algorithm}
\usepackage{algpseudocode}

\usepackage{hyperref}
\hypersetup{hidelinks}

\usepackage[capitalise,noabbrev]{cleveref}

\crefname{figure}{Fig.}{Figs.}
\Crefname{figure}{Fig.}{Figs.}

\theoremstyle{plain}
\newtheorem{theorem}{Theorem}[section]
\newtheorem{lemma}[theorem]{Lemma}
\newtheorem{corollary}[theorem]{Corollary}
\newtheorem{proposition}[theorem]{Proposition}
\newtheorem{conjecture}[theorem]{Conjecture}

\theoremstyle{definition}

\theoremstyle{remark}
\newtheorem{remark}[theorem]{Remark}

\crefname{conjecture}{Conjecture}{Conjectures}
\Crefname{conjecture}{Conjecture}{Conjectures}
\crefname{remark}{Remark}{Remarks}
\Crefname{remark}{Remark}{Remarks}

\newcommand{\Lam}{\Lambda}
\newcommand{\F}{\mathcal{F}}
\newcommand{\E}{\mathbb{E}}
\newcommand{\R}{\mathbb{R}}
\DeclareMathOperator{\sgn}{sgn}

\begin{document}

\title{Certified minimax normal-scale mixtures for the logistic-normal integral}
\author{Ruijie Li\thanks{Health Services and Outcomes Research, National Healthcare Group, 1
    Mandalay Road, Singapore 308205, Singapore.
    Emails: \href{mailto:ruijie.li@nhghealth.com.sg}{\nolinkurl{ruijie.li@nhghealth.com.sg}}
    and \href{mailto:ruijie@leefy.net}{\nolinkurl{ruijie@leefy.net}}.}}
\date{27 November 2026}

\maketitle

\begin{abstract}
The logistic distribution function is a scale mixture of normal
    distribution functions, and Monahan and Stefanski approximate it by the
    \(K\)-term mixture \(\Lambda(x) \approx \sum_{i} w_i \Phi(s_i x)\) chosen to
    minimise the uniform error on the real line. Each such mixture reduces the
    logistic-normal integral to a sum of \(K\) normal distribution functions
    with the mixture's uniform error bound. This integral is the expectation of
    a logistic sigmoid over a normal distribution in statistics and the
    Gauss--Fermi integral in semiconductor physics. Monahan and Stefanski
    tabulate constants for \(K \le 8\) and conjecture that the optimal error
    equioscillates at \(2K\) points, while above \(K = 8\) the Remez exchange
    and local solvers stall at nonequioscillating stationary points. We prove
    that the tangent space of the mixture family is a Chebyshev system, which
    forces \(2K\)-point equioscillation at every interior minimiser and makes
    any such solution the unique global minimax, and we prove interior
    attainment for \(K \le 11\). A certification algorithm audits a candidate
    mixture in multiple-precision arithmetic, using a proved bound on the
    number of extrema as its completeness test, and three search pipelines reach
    \(K = 9, 10, 11\), the last with certified worst-case error
    \(1.26 \times 10^{-11}\). A signed \(K = 22\) mixture carries a relaxed
    certificate at \(3.605 \times 10^{-15}\). In causal mediation with a binary
    outcome, the certified sum needs no rare-outcome assumption and replaces
    stochastic integration with one closed-form evaluation.
\end{abstract}



\section{Introduction}
\label{sec:intro}

The logistic distribution function \(\Lam(x) = 1/(1 + e^{-x})\) is an exact
continuous scale mixture of normal distribution
functions~\cite{andrews-scalemi-1974, stefansk-normals-1991},
\(\Lam(x) = \int_0^\infty \Phi(s x)\, q(s)\, ds\), with \(q\) the
Kolmogorov--Smirnov mixing density and \(\Phi\) the standard normal distribution
function. Monahan and Stefanski~\cite{monahan-normals-1992} discretised the
mixture to \(K\) atoms, which gives the approximation family
\begin{equation} \label{eq:mixture}
    f_K(x; w, s) = \sum_{i=1}^K w_i\, \Phi(s_i x),
    \qquad w_i \ge 0,\ \sum_i w_i = 1,\ s_i > 0,
\end{equation}
and chose its parameters to minimise the worst-case error over the real line,
\(\Delta_K^\star = \inf_{(w, s)} \sup_{x} |\Lam(x) - f_K(x; w, s)|\).

Monahan and Stefanski tabulated the minimax weights and scales for
\(K \le 8\) as fixed points of the Remez exchange
algorithm~\cite{remez-determin-1934}, with a worst-case error
\(\Delta_8^\star \approx 2.1 \times 10^{-9}\). They left two questions open. First,
their \emph{equioscillation conjecture}, that the error of the optimal \(K\)-term
mixture has exactly \(2K\) alternating extrema of equal magnitude on
\((0, \infty)\), was a numerical observation without proof. The classical
Chebyshev alternation theorem and linear Remez
theory~\cite{watson-approxim-2000} do not apply, because the mixture family is
nonlinear in the scales. Second, the tabulation did not reach beyond \(K = 8\).
Above \(K = 8\) the minimax error falls below \(10^{-9}\), and the
equioscillation system develops nonequioscillating stationary
points, algebraic saddles, and narrow anisotropic basins on which Remez and
Newton iterations stall or diverge.

The approximation was built for the \emph{logistic-normal integral}
\begin{equation} \label{eq:logistic-normal}
    I(\eta, v) = \int_{-\infty}^{\infty} \Lam(\eta + m)\,
        \frac{1}{\sqrt{2 \pi v}}\, \exp\!\left( -\frac{m^2}{2 v} \right) dm,
\end{equation}
the expectation of a logistic sigmoid over a normal distribution, which has no
elementary closed form. Since each normal term \(\Phi(s_i(\eta + m))\) integrates
in closed form against the Gaussian density, replacing \(\Lam\) by \(f_K\)
reduces \eqref{eq:logistic-normal} to the finite sum
\begin{equation} \label{eq:int-from-mixture}
    I(\eta, v) \approx \hat I(\eta, v) = \sum_{i=1}^K w_i\,
        \Phi\!\left( \frac{s_i \eta}{\sqrt{1 + s_i^2 v}} \right),
    \qquad |I - \hat I| \le \sup_x |\Lam(x) - f_K(x; w, s)|
\end{equation}
uniformly over \((\eta, v) \in \R \times (0, \infty)\). A minimax mixture
therefore turns \eqref{eq:logistic-normal} into a branch-free closed-form
evaluation carrying a parameter-free error bound.

The integral \eqref{eq:logistic-normal} recurs wherever a logistic link is
averaged over a Gaussian latent variable. Logistic regression with covariate
measurement error and with random effects drove the original
approximations~\cite{monahan-normals-1992}. The integral also enters
the simulation of short-rate and credit models with lognormally distributed
rates, where it solves the heat equation with a logistic initial
condition~\cite{pirjol-logistic-2013}. In the simulation of organic-semiconductor
devices it is the Gauss--Fermi integral, the charge-carrier density of a
Gaussian density of states~\cite{paasch-chargec-2010, farrell-comparis-2018}.
There \(\eta\) is the Fermi level measured from the centre of that density and
\(\sqrt{v}\) is its width, both in units of the thermal energy. That
literature evaluates it by an analytical approximation with separate expressions
in the nondegenerate and degenerate regions~\cite{paasch-chargec-2010}, by
series of error functions~\cite{mehmetog-analytic-2014, selvaggi-analytic-2018},
or by numerical quadrature~\cite{patriarc-highlya-2019}. Drift-diffusion
models of these devices use the integral through two ratios, the diffusion
enhancement \(I / \partial_\eta I\) and the degeneracy factor
\(e^{-\eta} I\)~\cite{farrell-comparis-2018}, so they need the small values of
the nondegenerate regime to relative accuracy. \Cref{rmk:tilt} therefore pairs
the mixture with a known exponential-tilt identity, which bounds the relative
error across all \(\eta\). In causal mediation analysis with a
binary outcome and a continuous mediator, \eqref{eq:logistic-normal} is the
computational kernel. There the \emph{natural indirect effect}, the change in
outcome probability transmitted through the mediator, is a difference of two
such integrals.

When individual values of \eqref{eq:logistic-normal} suffice, the evaluation
problem is solved. Pirjol~\cite{pirjol-logistic-2013} derives an exact recursion
in \(\eta\) together with a convergent Poisson-summation series that evaluates
the integral at a given point \((\eta, v)\) to the accuracy of the working
arithmetic. In simulation, optimisation, and uncertainty propagation, however,
the integral is consumed as a function across the parameter space rather than at
a single point. Uncertainty propagation, for example, repeats an evaluation
across thousands of \((\eta, v)\) draws, and any statement about the resulting
interval inherits the error of the inner evaluations. For a
truncated series the convergence rate and truncation error depend on the
evaluation point, so a uniform guarantee over a sampling loop must control the
truncation order at every draw. Fixed-node Gauss--Hermite quadrature loses
accuracy as \(v\) grows because the integrand's poles at \(\eta + m = \pm i\pi\)
approach the real axis at distance \(\pi/\sqrt{2v}\) in the quadrature variable.
A finite pole expansion tends to the same limit at both ends of \(\R\), whereas
\(\Lam\) tends to \(0\) and \(1\), so it cannot approximate \(\Lam\) uniformly on
\(\R\). A closed form whose worst-case error is a single
certified constant, fixed before any parameter is seen, removes that
per-draw control. The finite sum also has closed-form derivatives and remains
a finite sum under a further Gaussian convolution (\cref{sec:multi}). The
mixture does not improve on the digits of a
pointwise evaluator at a fixed point. What it contributes is the closed form,
the certified uniform bound over all \((\eta, v)\), and the exact finite-sum
expressions that follow from it. \Cref{sec:pointwise} develops this comparison
against Gauss--Hermite quadrature, Gauss--Fermi approximations, rational pole
expansions, and the Poisson-summation series.

This paper resolves both open questions of Monahan and Stefanski over the
certified range. The conjecture stood for
thirty-four years because the family lies outside the linear Chebyshev theory.
The Haar property of the mixture tangent space supplies the missing
equioscillation mechanism, and the mechanism extends to other atom families
with the same tangent structure (\cref{rmk:atoms}). Our contributions are as
follows.
\begin{enumerate}
    \item We prove the Monahan--Stefanski equioscillation conjecture at every
        interior minimiser (\cref{thm:haar,cor:equiosc}). The tangent space of the
        mixture manifold is a Chebyshev (Haar) system, which through a
        Kolmogorov-variational argument forces the minimax error to equioscillate
        at exactly \(2K\) points. The same zero-counting engine yields global
        uniqueness, making any equioscillating interior solution the unique global
        minimax (\cref{cor:global}), since the family is varisolvent in Rice's
        sense. Interior attainment is proved for
        \(K = 1, \ldots, 11\) from the certified constants (\cref{prop:interior}),
        so the conjecture holds unconditionally over the certified range, with
        full proofs in \cref{sec:S-setup} to~\cref{sec:S-symmetry}.
    \item We establish a tangent-Wronskian identity (\cref{thm:wronskian}), proved
        symbolically for \(K \le 5\), verified exactly through \(K = 12\), and
        confirmed by direct Sturm evaluation at every certified solution. The
        identity supplies an explicit local-uniqueness certificate independent of
        the global argument and guides the scale-separation penalty of the
        \(K = 11\) search.
    \item We state the certification algorithm (\cref{alg:certify}), whose scan
        for extrema is complete by the extremum bound of \cref{thm:S-upper}, and
        with it certify the minimax constants \((w^\star, s^\star, \Delta_K^\star)\)
        for \(K = 9, 10, 11\) (\cref{sec:certification}). Each case requires a
        different global-search pipeline: covariance matrix adaptation evolution
        strategy (CMA-ES)~\cite{hansen-complete-2001} with Remez and
        multiple-precision Newton polish at \(K = 9\), perturbed-reference
        iterated Remez at \(K = 10\), and
        Sturm-seeded CMA-ES with chain Levenberg--Marquardt
        (LM)~\cite{levenber-methods-1944, marquard-algorith-1963} and a final
        100-digit polish at \(K = 11\). The certified worst-case
        error at \(K = 11\) is \(\Delta_{11}^\star = 1.258 \times 10^{-11}\).
    \item We map the structure beyond \(K = 11\) (\cref{sec:beyond}): a conjectured
        plateau at \(K = 12, 13\), certified mixtures refuting its extension to
        \(K = 14, 15\), relaxed
        Bernstein--Akhiezer~\cite{akhiezer-theorya-1956} certificates through
        \(K = 17\), a proof that \(\Delta_K^\star \to 0\) as \(K \to \infty\)
        (\cref{prop:limit}), and a signed-class \(K = 22\) frontier that carries a
        certified strict signed \(30\)-point equioscillation.
    \item We demonstrate the consequence for the motivating application
        (\cref{sec:mediation}). In causal mediation the certified mixture needs none
        of the rarity that the closed forms of Valeri and VanderWeele assume, and it
        replaces the stochastic simulation or fixed-grid quadrature used elsewhere
        with a single certified closed-form sum, delivering asymmetric
        quasi-Bayesian intervals in under \(15\) milliseconds per \(10^4\) draws.
\end{enumerate}

\section{The certified minimax normal-scale mixture}
\label{sec:mixture}

Certifying a candidate mixture and promoting a solver fixed point to the global
minimax require three structural properties established in this section. First,
the tangent space of the mixture manifold is a Chebyshev system
(\cref{thm:haar}), which forces the minimax error to equioscillate at \(2K\)
points (\cref{cor:equiosc}). A mixed-Laplace variation-diminishing bound caps
the local extrema on \((0, \infty)\) of the error of any \(K\)-term mixture at
\(2K\) (\cref{thm:S-upper}).
Second, a zero-counting difference argument makes any \(2K\)-equioscillating
interior mixture the unique global minimax (\cref{cor:global}). Third, interior
attainment for \(K \le 11\) (\cref{prop:interior}) removes the interiority
hypothesis over the certified range. A separate tangent-Wronskian identity
(\cref{thm:wronskian}) supplies an explicit finite-\(K\) local-uniqueness
certificate and the scale-separation penalty used in the \(K = 11\) search.

\subsection{Equioscillation, attainment, and uniqueness}
Let \(\F_K \subset C^\infty(\mathbb{R})\) be the manifold of \(K\)-term
normal-scale mixtures \eqref{eq:mixture} with strictly positive weights and
ordered positive distinct scales \(0 < s_1 < \dots < s_K\), and write
\(\overline{\F_K}\) for its compactification with \(w_i \ge 0\) and
\(s_i \in [0, \infty]\). At any interior point \(f \in \F_K\) the tangent space is
\begin{equation} \label{eq:tangent}
    V = T_f \F_K = \operatorname{span}\!\left\{
        \Phi(s_1 z), \ldots, \Phi(s_K z),\
        z\, \phi(s_1 z), \ldots, z\, \phi(s_K z)
    \right\},
\end{equation}
of dimension \(2K - 1\) after enforcing the linear constraint \(\sum_i a_i = 0\)
on the \(\Phi\)-coefficients induced by \(\sum_i w_i = 1\). Here \(\phi = \Phi'\)
is the standard normal density. Since \(\Lam(-z) = 1 - \Lam(z)\) and
\(\sum_i w_i = 1\), the error \(e(z) = \Lam(z) - f(z)\) is odd on
\(\mathbb{R}\) and vanishes at the origin, so the minimax problem on
\(\mathbb{R}\) reduces to \((0, \infty)\). The tangent space collects the
first-order
changes of the mixture under small moves of its weights and scales, and so
governs whether the worst-case error can be reduced locally. In linear Chebyshev
approximation the Haar property of an \(n\)-dimensional approximating space
guarantees at least \(n + 1\) alternations at the optimum. For the nonlinear
family \(\F_K\) the tangent space assumes that role at each point, which is
why its Haar property is the structural fact on which the results below rest.

Every result of this section reduces to counting zeros, and the counts rest on
one classical bound.

\begin{lemma}[Zeros of generalised exponential sums] \label{lem:exp-sum-zeros}
    A generalised exponential sum \(\sum_{i=1}^m P_i(x)\, e^{\lambda_i x}\) with
    distinct real \(\lambda_i\) and \(\deg P_i = r_i\) has at most
    \(\bigl(\sum_i (r_i + 1)\bigr) - 1\) real zeros, counted with multiplicity. In
    particular \(\sum_{i=1}^N c_i e^{\mu_i x}\) has at most \(N - 1\), and
    \(\sum_{i=1}^K (c_i + d_i x)\, e^{\lambda_i x}\) at most \(2K - 1\).
\end{lemma}

The bound is classical~\cite{polya-ueberan-1913, schoenbe-polyafr-1951,
karlin-tchebych-1966}. Multiplying by \(e^{-\lambda_1 x}\), differentiating
\(r_1 + 1\) times to remove the first block, and applying Rolle's theorem gives
it by induction on \(\sum_i (r_i + 1)\).

\begin{theorem}[Haar property of the tangent space] \label{thm:haar}
    Let \(V\) be as in \eqref{eq:tangent}. Every nonzero \(h \in V\) has at most
    \(2K - 2\) distinct zeros in \((0, \infty)\). Equivalently, \(V\) is a
    Chebyshev (Haar) system of dimension \(2K - 1\) on \((0, \infty)\).
\end{theorem}

\begin{proof}
    Let \(h \in V\) be nonzero, and write \(h(z) = \sum_{i=1}^K a_i \Phi(s_i z) +
    \sum_{i=1}^K b_i\, z\, \phi(s_i z)\) with \(\sum_i a_i = 0\). Since
    \(\Phi(0) = 1/2\) and \(0 \cdot \phi(0) = 0\), \(h(0) = 0\), and since
    \(\Phi(s_i z) \to 1\) and \(z \phi(s_i z) \to 0\) as \(z \to \infty\),
    \(\lim_{z \to \infty} h(z) = 0\). Suppose \(h\) has \(N \ge 2K - 1\) distinct
    zeros in \((0, \infty)\). Adjoining boundary zeros at \(0\) and \(\infty\),
    Rolle's theorem gives \(h'\) at least \(N + 1\) zeros in \((0, \infty)\). After
    the substitution \(x = z^2\), \(\lambda_i = s_i^2 / 2\), zeros of \(h'(z)\) on
    \((0, \infty)\) biject with zeros of \(g(x) = \sum_i c_i\, \exp(-\lambda_i x) +
    \sum_i d_i\, x\, \exp(-\lambda_i x)\) on \((0, \infty)\), with \((c, d) \ne 0\)
    since the coefficient map is invertible for positive scales. By
    \cref{lem:exp-sum-zeros} this \(g\) has at most \(2K - 1\) zeros in
    \((0, \infty)\). Equivalently, the family \(\{\exp(-\lambda_i x),\, x\,
    \exp(-\lambda_i x)\}_{i=1}^K\) is a \(2K\)-dimensional Extended Complete
    Chebyshev system on \([0, \infty)\)~\cite{karlin-totalpo-1968,
    braess-nonlinea-1986}. Hence \(N + 1 \le 2K - 1\) and \(N \le 2K - 2\), a
    contradiction.
\end{proof}

\begin{corollary}[Equioscillation] \label{cor:equiosc}
    Let \(f^\star \in \F_K\) achieve \(\Delta_K^\star =
    \sup_{z \ge 0} |\Lam(z) - f^\star(z)|\) at an interior point. Then the error
    \(e^\star = \Lam - f^\star\) admits exactly \(2K\) alternating extrema in
    \((0, \infty)\), with \(|e^\star(z_j)| = \Delta_K^\star\) for
    \(j = 1, \ldots, 2K\) and \(\sgn\, e^\star(z_j) = (-1)^{j + 1}\, \sgn\,
    e^\star(z_1)\). By symmetry of \(\Lam\) and \(\F_K\) about \(z = 0\), the
    extended error on \(\mathbb{R}\) has exactly \(4K\) such extrema.
\end{corollary}

\begin{proof}[Proof sketch]
    The lower bound on the number of alternating extrema follows from
    \cref{thm:haar} via the Kolmogorov criterion for nonlinear Chebyshev
    approximation~\cite{cheney-introduc-1982, rice-approxim-1969}. If \(f^\star\)
    has \(N \le 2K - 1\) alternating extrema, the Haar property of \(V\) yields a
    tangent direction \(h \in V\) matching the sign of the error at every
    extremum. A parameter move realising \(h\) to first order, with a uniform
    second-order remainder, then decreases the error, contradicting minimaxity.
    The matching upper bound of \(2K\) holds for any \(K\)-term mixture with
    distinct positive scales (\cref{thm:S-upper}). The derivative \(e'\) of its
    error is a mixed Laplace transform in \(u = z^2/2\), consisting of a positive mixing
    density plus \(K\) atoms. A double-peel Rolle induction (the zero-counting
    form of the variation-diminishing property of the sign-regular Laplace
    kernel~\cite{karlin-totalpo-1968}) allows such a transform at most \(2K\)
    zeros, and every interior extremum of \(e^\star\) is a zero of
    \((e^\star)'\). The complete proof chain (setup, existence,
    the Haar property, the lower bound, the upper bound, uniqueness, and symmetry,
    \cref{sec:S-setup} to~\cref{sec:S-symmetry}, with an obstructed alternative
    total-positivity route in \cref{sec:S-alt-tp} and a proof-dependency index in
    \cref{sec:S-dag}) is given in the appendices.
\end{proof}

\Cref{cor:equiosc} is stated at an interior minimiser. The next result, proved in
\cref{sec:S-existence}, removes that hypothesis over the certified range. Attainment is
constructive, and no minimiser degenerates.

\begin{proposition}[Interior attainment for \(K \le 11\)] \label{prop:interior}
    For \(K = 1, \ldots, 11\) the minimax \(\Delta_K^\star\) is attained, and every
    minimiser is an interior point of \(\F_K\), with no vanishing weight and no
    two scales coinciding or escaping to \(0\) or \(\infty\).
\end{proposition}

\begin{proof}[Proof sketch]
    Attainment on the compactification \(\overline{\F_K}\) follows from compactness of
    the parameter domain and lower semicontinuity of the supremum norm under pointwise
    limits. A minimiser can leave the interior by two boundary mechanisms, and an
    induction on \(K\) from the \(K = 1\) base, which supplies the closure bound
    at level \(K - 1\), excludes both.

    First, a weight may vanish or two scales may coincide. For \(K \ge 2\) such a
    count-reducing minimiser would lie in \(\overline{\F_{K-1}}\) and incur an
    error of at least \(\Delta_{K-1}^\star > \Delta_K^\star\). That strict
    inequality requires no attainment assumption, since the exhibited
    \(K\)-mixture of \cref{sec:certification} bounds \(\Delta_K^\star\) from above
    and the audited alternation of the \((K-1)\)-mixture bounds
    \(\Delta_{K-1}^\star\) from below through the de la Vall\'ee Poussin inequality
    (\cref{prop:S-dlvp}). At \(K = 1\) no such minimiser exists, since
    \(\sum_i w_i = 1\) forces the single weight to one.

    Second, mass may escape to the boundary scales \(s = 0\) and \(s = \infty\),
    where it contributes a constant, or a step at the origin, and forces a supremum
    error of half that mass. For \(K \ge 2\) a comparison against
    \(\Delta_{K-1}^\star\) excludes this whenever
    \(\Delta_K^\star / \Delta_{K-1}^\star < (1 - 4 \Delta_K^\star)/2\). The
    certified ratios do not exceed \(0.19\), while the threshold exceeds \(0.498\)
    throughout. At \(K = 1\) the only boundary points are the constant
    \(\tfrac{1}{2}\) and the unit step, each with error \(\tfrac{1}{2}\).
\end{proof}

\begin{corollary}[Global uniqueness] \label{cor:global}
    For \(K = 1, \ldots, 11\) the best uniform approximation to \(\Lam\) from
    \(\F_K\) is unique, and any interior candidate whose error equioscillates at
    \(2K\) points is that unique global minimax.
\end{corollary}

\begin{proof}[Proof sketch]
    Two distinct mixtures in \(\F_K\) differ by a signed combination of at most
    \(2K\) distinct \(\Phi\)-atoms with zero coefficient sum. Such a combination
    has at most \(2K - 2\) zeros on \((0, \infty)\). Alternation of the error of
    one of the two at \(2K\) points at the minimax level forces at least
    \(2K - 1\) zeros, the
    equality case of the de la Vall\'ee Poussin argument (\cref{sec:S-uniqueness}).
    Classically, the same conclusion follows from
    varisolvency~\cite[Theorem~7-4]{rice-approxim-1969},
    \cite[III.4.6]{braess-nonlinea-1986}.
    Combining with \cref{prop:interior,cor:equiosc} gives the claim.
\end{proof}

\begin{figure}[htbp]
    \centering
    \begin{adjustbox}{max width=\textwidth}
    \begin{tikzpicture}
        \begin{axis}[
            width=14cm, height=5.4cm,
            axis lines=middle,
            enlarge x limits=false,
            xmin=0, xmax=13, ymin=-6.8, ymax=6.8,
            xlabel={\(z\)}, ylabel={\(e^\star(z)\)},
            xlabel style={at={(axis description cs:1,0.5)}, anchor=west},
            ylabel style={at={(axis description cs:0,1)}, anchor=south east},
            xtick={2,4,6,8,10,12},
            ytick={-5.9529,0,5.9529},
            yticklabels={\(-\Delta_5^\star\),\(0\),\(+\Delta_5^\star\)},
            tick label style={font=\footnotesize},
        ]
            \addplot[dashed, gray, line width=1pt, samples=2, domain=0:13] {5.9529};
            \addplot[dashed, gray, line width=1pt, samples=2, domain=0:13] {-5.9529};
            \addplot[black, very thick] table[x=z, y=e, col sep=comma] {fig-equiosc.dat};
            \addplot[only marks, mark=*, mark size=1.5pt, black]
                table[x=z, y=e, col sep=comma] {fig-equiosc-ext.dat};
        \end{axis}
    \end{tikzpicture}
    \end{adjustbox}
    \caption{The certified \(K = 5\) minimax error \(e^\star(z) = \Lam(z) -
        \sum_{i=1}^{5} w_i^\star \Phi(s_i^\star z)\) on \((0, \infty)\), computed from the
        certified constants of \cref{tab:cert}. In agreement with \cref{cor:equiosc}, it
        attains \(\pm\Delta_5^\star\) (dashed, \(\Delta_5^\star = 5.953 \times 10^{-7}\)) at
        exactly \(2K = 10\) points with strictly alternating signs (marked), and vanishes at
        the origin by oddness. The extrema spread apart as \(z\) grows, the outermost near
        \(z = 12\)}
    \label{fig:equiosc}
\end{figure}
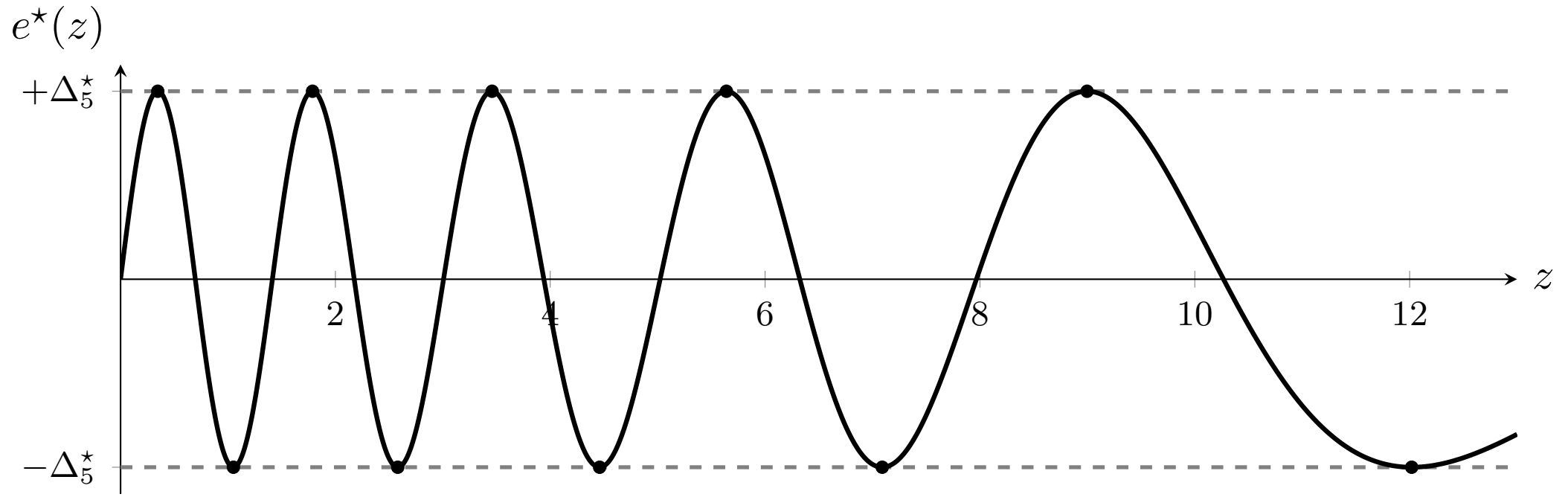

\Cref{cor:equiosc,prop:interior} together resolve the Monahan--Stefanski
conjecture unconditionally for \(K \le 11\), and \cref{cor:global} promotes every
certified equioscillating minimum of \cref{sec:certification} to the unique strict
global minimum, so any Remez-type fixed point gains a termination certificate.
\Cref{fig:equiosc} illustrates the resulting equioscillation.

\begin{remark}[Other atom families] \label{rmk:atoms}
    The proof of \cref{thm:haar} uses only that the derivative of the tangent
    space maps, under the substitution \(x = z^2\), onto an Extended Chebyshev
    space, a property not special to Gaussian atoms. It generalises to any
    smooth atom family \(\{G(s_i z)\}\) with even density \(G'\) for which the
    derivative block \(\{G'(s_i z),\, z\,G''(s_i z)\}\) transforms under
    \(x = z^2\) onto a known Extended Chebyshev space. That block collapses to
    the Gaussian form \(\{G'(s_i z),\, z^2 G'(s_i z)\}\) only when \(G''(u)\) is
    a constant multiple of \(u\,G'(u)\), which characterises the Gaussian, so
    any other atom is governed by the general block. Every \(G\) in that class
    inherits the Haar property of the tangent space, and with it the lower bound
    of at least \(2K\) alternations. The matching upper bound and the uniqueness
    argument require further structure, since the exact count rests on a continuous
    \(G\)-scale mixture with a sign-regular derivative kernel and uniqueness rests
    on a bound of \(2K - 2\) zeros for signed combinations of \(2K\) atoms. The
    logistic, with \(G'' = G'\,(1 - 2G)\), is a plausible member that the
    present proof does not settle.
\end{remark}

\subsection{Wronskian local uniqueness}
The uniqueness argument above works through the zero-counting structure of the
family as a whole. Independent verification also calls for an algebraic check
that can be evaluated directly at each certified scale vector \(s^\star\), and
whose scale factor measures proximity to scale coalescence during global search
(\cref{sec:pipelines}). The tangent-Wronskian identity supplies both. A cofactor
of the tangent Wronskian that vanishes nowhere on the ordered positive cone
certifies that the tangent basis is nondegenerate at every candidate,
independently of the equioscillation theorem. Let
\(Q(z; s)\) denote the \((2K-1) \times (2K-1)\) tangent-Wronskian cofactor of the
basis \(\{\Phi(s_i z), z\, \phi(s_i z)\}_{i=1}^K\) after deleting the order-\(0\)
derivative row and the largest-index column \(z\, \phi(s_K z)\). Rows are
ordered by ascending derivative order, and the surviving columns as the \(\Phi\)
block followed by the \(z\,\phi\) block, the ordering that fixes the sign of the
closed form in \cref{sec:S-uniqueness}.

\begin{theorem}[Structural pure-monomial identity, \(K \le 5\)] \label{thm:wronskian}
    For \(2 \le K \le 5\) and ordered positive distinct
    \(s_1 < s_2 < \dots < s_K\), the reduced tangent-Wronskian cofactor satisfies
    \(Q(z; s) = c_K(s)\, z^{(2K-1)(K-1)}\), where \(c_K(s)\), a signed product of a
    power of two, powers of the scales, and even powers of the pairwise differences
    of squared scales given in \cref{sec:S-uniqueness}, is nonzero on the ordered positive
    cone. The Wronskian is therefore nonsingular at every interior point of \(\F_K\).
\end{theorem}

The proof is a full symbolic computation. A Bareiss fraction-free
elimination~\cite{bareiss-sylveste-1968} for \(K \le 4\) and a staged suffix-minor expansion in exact
integer arithmetic for \(K = 5\) return the reduced determinant as an exact
multivariate polynomial, which factors as stated (\cref{sec:S-uniqueness}). That an
equioscillating candidate is a strict local minimum is a separate conclusion,
reached in \cref{sec:S-uniqueness} by the uniqueness corollary, which draws strictness from
the Haar and de la Vall\'ee Poussin chain rather than from \(Q\).

\begin{conjecture}[Pure-monomial identity for general \(K\)]
\label{conj:wronskian}
    The identity of \cref{thm:wronskian} holds for every \(K \ge 2\).
\end{conjecture}

\Cref{conj:wronskian} is verified exactly at the integer scales
\(s = (1, 2, \ldots, K)\) for \(6 \le K \le 12\). The Bareiss determinant and the
closed form evaluate to the same integer, which exceeds \(400\)
digits at \(K = 12\). The identity is a checkable certificate rather than a new
uniqueness route. For each certified
mixture of \cref{sec:certification} the nonvanishing of \(Q\) at the certified
\(s^\star\) is in addition confirmed by a direct Sturm real-root count carried
out in exact rational arithmetic, so the local-uniqueness certificate there does not rely on
\cref{conj:wronskian}.

\section{Certified solutions for \texorpdfstring{\(K = 1, \ldots, 11\)}{K = 1 to 11}}
\label{sec:certification}

The attainment and uniqueness results of \cref{sec:mixture} rest on exhibiting
mixtures whose uniform errors decrease strictly with \(K\) and equioscillate at
\(2K\) points. This section reports those certified constants for
\(K = 1, \ldots, 11\), states the certification and evaluation algorithms,
details the search pipelines that solved \(K = 9, 10, 11\), and records how
every result is reproduced.

\Cref{tab:cert} reports the certified minimax constants for
\(K = 1, \ldots, 11\). For \(K \le 8\) these are the Monahan--Stefanski
cases~\cite{monahan-normals-1992}, re-solved here in 100-digit arithmetic. At
\(K = 8\), the largest case, our constants agree with their published values to
eight significant figures. Evaluated as printed, the published constants give a
supremum that exceeds \(\Delta_8^\star\) by \(4.9 \times 10^{-7}\) in relative
terms, or by \(1.9 \times 10^{-7}\) once the \(10^{-15}\) by which their
weights overshoot unity is removed. For \(K = 9, 10, 11\) the constants
are new. Every entry satisfies \cref{cor:equiosc} in 100-significant-digit
arithmetic. All \(2K\) extrema have equal magnitude to relative spread below
\(10^{-50}\) for \(K = 9, 10, 11\), and the signs strictly alternate.

\begin{table}[htbp]
    \footnotesize
    \caption{Certified minimax constants \(\Delta_K^\star\) for
        \(K = 1, \ldots, 11\), with the search pipeline that closed each case.
        M\&S denotes the cases of Monahan and
        Stefanski~\cite{monahan-normals-1992}, re-solved under the 100-digit
        audit and agreeing with their published constants to eight significant
        figures at \(K = 8\), the largest case. CMA-ES denotes covariance matrix
        adaptation evolution strategy, LM denotes Levenberg--Marquardt, and MPFR
        denotes the GNU multiple-precision floating-point library. Displayed
        values are rounded to the digits shown. The audited
        full-precision weights and scales are in the accompanying repository, in
        \texttt{inst/certs/minimax\_k01.csv} through \texttt{minimax\_k11.csv},
        and the audited bounds themselves in \texttt{inst/certs/deltas.csv}.}
    \label{tab:cert}
    \centering
    \begin{tabular}{rlcl}
        \toprule
        \(K\) & Certified \(\Delta_K^\star\) (100-digit audit) & Extrema & Pipeline \\
        \midrule
        1  & \(9.4573 \times 10^{-3}\)                & 2  & M\&S                              \\
        2  & \(5.1013 \times 10^{-4}\)                & 4  & M\&S                              \\
        3  & \(4.3559 \times 10^{-5}\)                & 6  & M\&S                              \\
        4  & \(4.7167 \times 10^{-6}\)                & 8  & M\&S                              \\
        5  & \(5.9529 \times 10^{-7}\)                & 10 & M\&S                              \\
        6  & \(8.3872 \times 10^{-8}\)                & 12 & M\&S                              \\
        7  & \(1.2857 \times 10^{-8}\)                & 14 & M\&S                              \\
        8  & \(2.1086 \times 10^{-9}\)                & 16 & M\&S                              \\
        9  & \(3.656450583166549933 \times 10^{-10}\) & 18 & CMA-ES + Remez + MPFR Newton          \\
        10 & \(6.646923022720961975 \times 10^{-11}\) & 20 & Perturbed-reference Remez + MPFR Newton \\
        11 & \(1.258437696716086980 \times 10^{-11}\) & 22 & Sturm-seeded CMA-ES + chain LM + polish \\
        \bottomrule
    \end{tabular}
\end{table}

\begin{figure}[htbp]
    \centering
    \begin{adjustbox}{max width=\textwidth}
    \begin{tikzpicture}
        \begin{axis}[
            width=14cm, height=6.6cm,
            axis lines=left,
            ymode=log,
            xmin=0, xmax=24.4, ymin=1e-12, ymax=3e-2,
            xtick={0,4,8,12,16,20},
            ytick={1e-2,1e-4,1e-6,1e-8,1e-10,1e-12},
            xlabel={\(z\)}, ylabel={\(|e_K^\star(z)|\)},
            label style={font=\footnotesize},
            tick label style={font=\footnotesize},
            clip mode=individual,
        ]
            \pgfplotsforeachungrouped \Kval/\Dval in {1/9.4573e-3, 2/5.1013e-4,
                3/4.3559e-5, 4/4.7167e-6, 5/5.9529e-7, 6/8.3872e-8, 7/1.2857e-8,
                8/2.1086e-9, 9/3.6565e-10, 10/6.6469e-11, 11/1.2584e-11} {
                \addplot[dashed, gray, line width=0.5pt, samples=2, domain=0:22] {\Dval};
                \addplot[black, line width=0.35pt]
                    table[x=z, y=a\Kval] {fig-convergence.dat};
            }
            \node[font=\tiny, anchor=west] at (axis cs: 22.4, 9.4573e-3)  {\(K{=}1\)};
            \node[font=\tiny, anchor=west] at (axis cs: 22.4, 5.1013e-4)  {\(K{=}2\)};
            \node[font=\tiny, anchor=west] at (axis cs: 22.4, 4.3559e-5)  {\(K{=}3\)};
            \node[font=\tiny, anchor=west] at (axis cs: 22.4, 4.7167e-6)  {\(K{=}4\)};
            \node[font=\tiny, anchor=west] at (axis cs: 22.4, 5.9529e-7)  {\(K{=}5\)};
            \node[font=\tiny, anchor=west] at (axis cs: 22.4, 8.3872e-8)  {\(K{=}6\)};
            \node[font=\tiny, anchor=west] at (axis cs: 22.4, 1.2857e-8)  {\(K{=}7\)};
            \node[font=\tiny, anchor=west] at (axis cs: 22.4, 2.1086e-9)  {\(K{=}8\)};
            \node[font=\tiny, anchor=west] at (axis cs: 22.4, 3.6565e-10) {\(K{=}9\)};
            \node[font=\tiny, anchor=west] at (axis cs: 22.4, 6.6469e-11) {\(K{=}10\)};
            \node[font=\tiny, anchor=west] at (axis cs: 22.4, 1.2584e-11) {\(K{=}11\)};
        \end{axis}
    \end{tikzpicture}
    \end{adjustbox}
    \caption{The certified minimax error magnitudes
        \(|e_K^\star(z)| = |\Lam(z) - \sum_{i=1}^{K} w_i^\star \Phi(s_i^\star z)|\)
        for the whole certified range \(K = 1, \ldots, 11\), computed in
        \(256\)-bit arithmetic from the constants of \cref{tab:cert}. Grey dashed
        lines mark the certified \(\Delta_K^\star\). In agreement with
        \cref{cor:equiosc}, curve \(K\) touches its envelope at exactly \(2K\)
        points, and every one of them lies inside the window, with the outermost
        moving outward from \(z = 3.48\) at \(K = 1\) to \(z = 21.50\) at
        \(K = 11\). The downward spikes are zero crossings of \(e_K^\star\),
        clipped by the log axis. Past its outermost extremum each curve settles
        onto the common asymptote \(1 - \Lam(z) \sim e^{-z}\), since the Gaussian
        tails have by then decayed below the logistic one, so the curves join that
        asymptote in order, \(K = 1\) first and \(K = 11\) last. The envelopes
        fall geometrically at the ratios given in the text, so the eleven
        certificates span nearly nine orders of magnitude}
    \label{fig:convergence}
\end{figure}
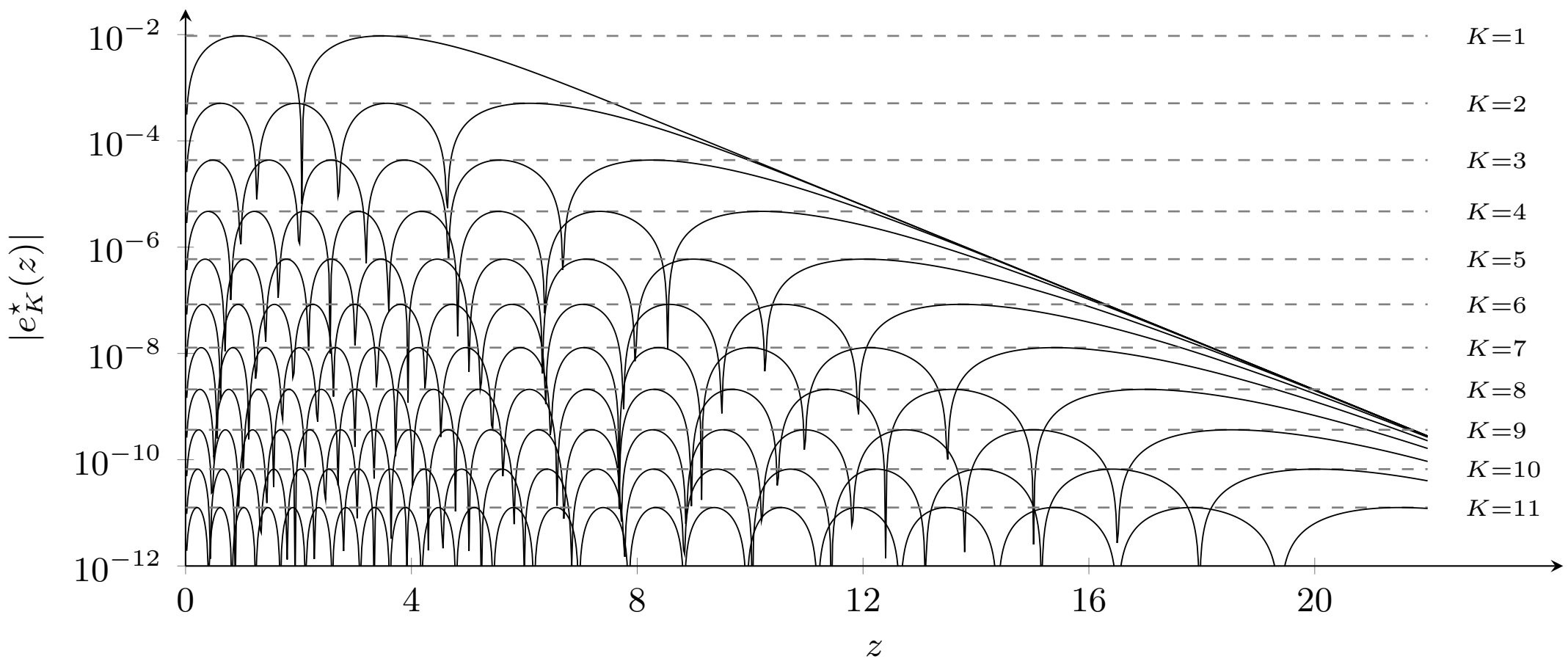

The ratio \(\Delta_{K+1}^\star / \Delta_K^\star\) takes the values \(0.054\),
\(0.085\), \(0.108\), \(0.126\), \(0.141\), \(0.153\), \(0.164\), \(0.173\),
\(0.182\), and \(0.189\) for \(K = 1 \to 2, \ldots, 10 \to 11\). The sequence
increases monotonically towards unity but remains strictly below it through
\(K = 11\). A least-squares fit of \(\log \Delta_K^\star\) on
\(\{1, K, K^2, \log K, K \log K\}\) over \(K = 1, \ldots, 10\), extrapolated one
step, predicts \(\Delta_{11}^\star \approx 1.256 \times 10^{-11}\), which the
certified value matches to within \(0.19\%\), an independent consistency check
on the certification. \Cref{fig:convergence} plots all eleven certified error
curves against the envelopes those constants define.

\begin{remark}[What the certificate guarantees]
    The certified value \(\Delta_K^\star\) is the supremum of
    \(|\Lam(z) - \sum_i w_i \Phi(s_i z)|\) audited in 100-digit arithmetic for
    the specific weights and scales shipped, an upper bound the deployed mixture
    is verified to meet. With \(\hat I\) as in \eqref{eq:int-from-mixture}, the
    accuracy guarantee \(|I - \hat I| \le \Delta_K^\star\) holds for that
    mixture whether or not it is the global minimax. Global optimality, that no
    \(K\)-term mixture attains a smaller uniform error, is the separate and
    stronger statement proved for \(K \le 11\) by
    \cref{cor:equiosc,cor:global}. The deployed error bound does not depend on
    it.

    Conversely, the de la Vall\'ee Poussin inequality, proved
    for this family by an elementary difference argument (\cref{sec:S-uniqueness}), turns the
    verified alternation at levels within relative spread \(\varepsilon\) of the
    audited value into the lower bound \((1 - \varepsilon)\) times that value.
    Each certified constant therefore gives a two-sided enclosure of the true
    minimax at the stated spread.
\end{remark}

\subsection{The certification and evaluation algorithms}
\label{sec:algorithms}
A search returns a candidate, and a candidate becomes a certificate only when the
supremum of its error has been computed and its alternation verified.
\Cref{alg:certify} is that step. It is the same for every \(K\) and for every
search pipeline, and it takes the candidate's weights and scales as exact
decimals, so the bound it returns belongs to the constants as distributed.

Since the distributed weights are stored as finite decimals, write
\(d = 1 - \sum_i w_i\) for their unit-sum defect, which is zero for exact
constants and below \(10^{-39}\) in magnitude for every set distributed here.
Then \(e(-z) = d - e(z)\), so \(e'\) is even, the critical points on the
negative half-line mirror those on the positive one, and \(e(z)\) tends to
\(0\) as \(z \to -\infty\) and to \(d\) as \(z \to +\infty\). The supremum of
\(|e|\) over the real line is therefore attained at a critical point on
\((0, \infty)\) or at its mirror image on \((-\infty, 0)\), unless it is at
most \(|d|\).

To locate all critical points on \((0, \infty)\), \cref{alg:certify} scans
\(e'\) for sign changes on a finite grid \(G \subset (0, Z]\). For a general
function such a scan is incomplete, since a pair of critical points falling
between two neighbouring grid nodes or beyond \(Z\) would go undetected and
understate the supremum. \Cref{thm:S-upper} removes that obstacle for this
family. The error of any \(K\)-term mixture with distinct positive scales has
at most \(2K\) local extrema on \((0, \infty)\), so a grid scan that brackets
\(2K\) sign changes of \(e'\) has provably found every local extremum.

\begin{algorithm}[htbp]
    \caption{Certification of a \(K\)-term mixture}
    \label{alg:certify}
    \begin{algorithmic}[1]
        \Require real weights \(w_i\) and distinct scales \(s_i > 0\), as exact
            decimals with \(d = 1 - \sum_i w_i\), a grid \(G \subset (0, Z]\), and
            two working precisions \(p_{\mathrm{scan}} < p\) in significant digits
        \Ensure a bound \(\Delta \ge \sup_{z \in \R} |e(z)|\) for
            \(e(z) = \Lam(z) - \sum_i w_i \Phi(s_i z)\), and the spread
            \(\varepsilon\) of its alternation
        \State Evaluate \(e'(z) = \Lam'(z) - \sum_i w_i s_i \phi(s_i z)\) on \(G\)
            at precision \(p_{\mathrm{scan}}\), and keep each pair of neighbouring
            grid points at which \(e'\) takes opposite signs
            \Comment{brackets}
        \State Bisect each bracket at precision \(p\) to a critical point
            \Comment{\(z_1 < \cdots < z_n\)}
        \State Require \(n = 2K\) and
            \(\sgn e(z_{j+1}) = -\sgn e(z_j)\) for \(j = 1, \ldots, 2K - 1\)
            \Comment{complete by \cref{thm:S-upper}}
        \State \(\Delta \gets \max_j \max\{|e(z_j)|, |d - e(z_j)|, |d|\}\), rounded
            up at the last digit reported
            \Comment{both half-lines}
        \State \(\varepsilon \gets 1 - \min_j |e(z_j)| / \Delta\)
            \Comment{\((1 - \varepsilon) \Delta - |d| \le \Delta_K^\star \le
            \Delta + |d|\)}
    \end{algorithmic}
\end{algorithm}

Stage 1 brackets the zeros of \(e'\), stage 2 refines the critical points,
stage 3 checks completeness and alternation, stage 4 bounds \(|e|\) on both
half-lines, and stage 5 forms the two-sided enclosure for positive weights.
Dividing the weights by their sum gives a member of \(\F_K\) whose error
differs from \(e\) by at most \(|d|\). The upper side of the enclosure holds
since that member belongs to \(\F_K\), and the lower side is the de la Vall\'ee Poussin
inequality of \cref{prop:S-dlvp} applied to its \(2K\) alternating values. The
bounds distributed with the constants were recomputed by \cref{alg:certify}
with \(G\) the union of \(16{,}000\) equally spaced points on \((0, 80]\) and
\(3{,}001\) logarithmically spaced points on \([10^{-4}, 1]\), at
\(p_{\mathrm{scan}} = 40\) and \(p = 140\). The arithmetic is
multiple-precision floating point rather than interval arithmetic, so the
certificate is a verification at the stated precision
(\cref{sec:reproducibility}). The relaxed certificates of \cref{sec:beyond}
have fewer than \(2K\) critical points, namely \(22\), \(22\), \(21\), \(21\),
and \(31\) at \(K = 14, 15, 16, 17, 22\). For those mixtures stages 3 and 5 do
not apply, and the completeness of the scan rests on the density of \(G\).

\Cref{alg:evaluate} is the use the certificate is built for. It restates
\eqref{eq:int-from-mixture} as a procedure.

\begin{algorithm}[htbp]
    \caption{Evaluation of the logistic-normal integral from a certified mixture}
    \label{alg:evaluate}
    \begin{algorithmic}[1]
        \Require \(\eta \in \R\), \(v \ge 0\), and a certified mixture
            \((w, s)\) with bound \(\Delta\)
        \Ensure \(\hat I\) with \(|I(\eta, v) - \hat I| \le \Delta\)
        \For{\(i = 1, \ldots, K\)}
            \State \(t_i \gets s_i \eta / \sqrt{1 + s_i^2 v}\)
        \EndFor
        \State \(\hat I \gets \sum_{i=1}^{K} w_i \Phi(t_i)\)
    \end{algorithmic}
\end{algorithm}

The evaluation costs \(K\) values of \(\Phi\). It has no branch, no iteration,
and no tuning parameter, and its error bound is the same constant at every
\((\eta, v)\), which is what the sampling loops of \cref{sec:mediation} depend on.

\begin{remark}[Relative accuracy at small values] \label{rmk:tilt}
    The bound of \cref{alg:evaluate} is absolute. As \(\eta \to -\infty\) the
    integral decays as \(e^{\eta + v/2}\) while \(\hat I\) decays at a Gaussian
    rate, so the relative error of \(\hat I\) tends to one there. An application
    that needs small values to a fixed relative
    accuracy, such as the carrier density of a nondegenerate semiconductor, is
    therefore not served by \(\Delta\) alone. A known exact identity supplies
    what is missing. Since \(\Lam(x) = e^{x} \Lam(-x)\), the exponential tilt
    \(E[e^{Y} g(Y)] = e^{\eta + v/2} E[g(Y + v)]\) of the normal argument
    \(Y \sim N(\eta, v)\) gives
    \begin{equation} \label{eq:tilt}
        I(\eta, v) = e^{\eta + v/2}\, I(-\eta - v, v).
    \end{equation}
    Crouch and Spiegelman~\cite{crouch-evaluati-1990} state it among their
    symmetry relations, and it is the recursion of
    Pirjol~\cite{pirjol-logistic-2013} combined with the reflection
    \(I(-\eta, v) = 1 - I(\eta, v)\). Its instance at \(\eta = -v\), the exact
    value \(I(-v, v) = e^{-v/2} / 2\), marks the passage from the nondegenerate
    regime to the degenerate one in the semiconductor
    literature~\cite{paasch-chargec-2010}. Gory, Craigmile, and
    MacEachern~\cite{gory-classge-2021} apply the recursion before the
    eight-term mixture and assess the gain in accuracy numerically. With the
    certified \(\Delta\), one application yields a bound on the relative error.
    For \(\eta < -v/2\), evaluate the right-hand side by \cref{alg:evaluate}. The
    error of the product is at most \(e^{\eta + v/2} \Delta\). The argument
    \(-\eta - v\) exceeds \(-v/2\) and \(I\) increases in its first argument, so
    \(I(\eta, v) \ge e^{\eta + v/2} I(-v/2, v)\), and the relative error is at
    most \(\Delta / I(-v/2, v)\). For \(\eta \ge -v/2\) the plain evaluation
    meets the same bound, which therefore holds at every \(\eta\) and depends on
    \(v\) alone. At \(K = 11\) it is \(3.2 \times 10^{-11}\) at \(v = 1\) and
    \(6.1 \times 10^{-9}\) at \(v = 36\), since \(I(-v/2, v)\) falls from
    \(0.40\) to \(2.1 \times 10^{-3}\) between those variances, so the guarantee
    weakens as \(v\) grows. The plain sum of \cref{alg:evaluate} remains the default. For
    weights summing to one it satisfies
    \(\hat I(\eta, v) + \hat I(-\eta, v) = 1\) identically, which the tilted
    evaluation forgoes, and the effect expressions of \cref{sec:mediation} are
    stated for it.
\end{remark}

\subsection{Pipelines and failure modes}
\label{sec:pipelines}
The three new cases share one difficulty and needed three different remedies.
At candidate extrema \(0 < z_1 < \dots < z_{2K}\) with
\(\sigma = \sgn e(z_1)\), the equioscillation system imposes
\(e(z_j; w, s) = (-1)^{j+1}\sigma\Delta\) for \(j = 1, \ldots, 2K\) under
\(\sum_i w_i = 1\), coupled with \(e'(z_j; w, s) = 0\) in the full-system
polish. Write \(F\) for its residual in the weights, the scales, and \(\Delta\),
and \(J\) for the Jacobian of \(F\). At every \(K \ge 9\) the squared residual
\(\|F\|^2\) has saddles and nonequioscillating local minima that Newton,
Levenberg--Marquardt, and Nelder--Mead reach readily and leave only with
difficulty. The local geometry gives no first-order signal towards the
equioscillating global minimum that \cref{cor:equiosc} guarantees. The remedy
at each \(K\) injects a nonlocal move matched to that \(K\)'s basin geometry
(\cref{tab:remedies}).

\paragraph{\texorpdfstring{\(K = 9\)}{K = 9}: Global escape} Standard local and
global searches all converged on a nonequioscillating stationary point at
\(\Delta \approx 6.72 \times 10^{-10}\). These searches included Newton,
Nelder--Mead~\cite{nelder-simplex-1965}, differential evolution, Sobol and
extreme-range seeding, random basin-hopping over \(10^6\) samples, and
enumeration of \(2 \times 10^6\) nine-scale subsets from a log-spaced
dictionary. At that point, the error reaches \(99\%\) of its supremum at only
\(12\) of the \(18\) extrema that equioscillation requires. Every failed
variant's covariance or step size was isotropic in log-scale space, while the
target basin is strongly anisotropic (certified scale ratio
\(s_{\max}^\star / s_{\min}^\star \approx 6.3\)). Running covariance matrix
adaptation evolution strategy (CMA-ES)~\cite{hansen-complete-2001} on the
\(9\)-dimensional log-scale space with \(50\) restarts located the correct
basin, since its adapted covariance breaks that symmetry. A Remez warm start at
double precision reduced the error to \(\Delta \approx 3.67 \times 10^{-10}\)
at relative spread \(1.1\%\). A \(256\)-bit Newton polish in MPFR then
converged in \(13\) iterations, and a \(115\)-digit full-system Newton solve
produced the 100-digit constants of \cref{tab:cert}.

\paragraph{\texorpdfstring{\(K = 10\)}{K = 10}: Spread reduction} The \(K = 9\)
pipeline reached the correct basin at \(\Delta = 6.7095 \times 10^{-11}\) with
\(20\) alternating extrema but spread \(2.7\%\), where Newton diverged and
Levenberg--Marquardt stalled with \(J^{\!\top} F = 0\), the signature
of an algebraic saddle. A wider CMA-ES sweep returned a marginally different
stationary point of the same structure at \(6.6952 \times 10^{-11}\), so no
better basin remained to find. The perturbed-reference iterated Remez strategy
treats the extrema as control parameters. It perturbs each \(z_j^\star \to
z_j^\star (1 + \delta)\) for \(\delta \in \pm\{1, 2, 3, 5, 8, 10, 20, 30\}\%\),
jointly perturbs pairs among the six extrema of smallest magnitude,
and solves the Newton equioscillation system with \(z_j\) fixed. Three rounds drove the spread from
\(2.7\%\) through \(0.33\%\) to \(0.016\%\), after which MPFR Newton and the
\(115\)-digit full-system Newton solve converged quadratically.

\paragraph{\texorpdfstring{\(K = 11\)}{K = 11}: Chain Levenberg--Marquardt}
The \(K = 11\) problem combined both failure modes. The minimax basin width in
log-scale is below \(5 \times 10^{-5}\), and a single full Newton step from the
\(22\)-extremum candidate described below exits it, \(\|F\|\) rising from
\(1.4 \times 10^{-12}\) to \(1.7 \times 10^{-7}\). The smallest minimax weight
falls from \(w_{\min}^{(10)\star} \approx 2.45 \times 10^{-4}\) to
\(w_{\min}^{(11)\star} \approx 1.03 \times 10^{-4}\), ruling out weight floors
calibrated on the preceding solutions. Twenty-nine strategies failed before
the successful pipeline. They fall into three classes, namely
alternation-penalised global searches, local Newton polishes in coordinate and
log-scale space, and \(K = 10 + 1\) atom insertions. The first class reached
\(22\) alternations only at order-one relative spread or not at all. The one
exception was a saddle of \(\|F\|^2\) at \(\Delta = 1.303 \times 10^{-11}\)
(spread \(7.4\%\)), from which the second class raised the spread to \(13.5\%\)
and \(11.7\%\). Every trial in the third class collapsed to a \(K = 10\)
attractor.

Two ingredients completed the certification. The first is a Sturm-seeded CMA-ES
run whose objective penalises a trial when \(\log |c_K(s)|\), the
tangent-Wronskian coefficient of \cref{conj:wronskian}, falls below a target,
steering the search towards well-separated scales. Re-auditing that run in
100-digit arithmetic revealed \(22\) extrema where the linear-programming grid
verifier, a double-precision routine that computes best weights at fixed scales
on a grid, had reported \(21\). The grid had merged an extremum of small
magnitude with its neighbour. The second is chain Levenberg--Marquardt with
adaptive damping \(\lambda\) in 100-digit arithmetic. It threaded the narrow
basin without exiting it across five segments of \(500\), \(500\), \(500\),
\(10^4\), and up to \(2 \times 10^5\) iterations, stopping once
\(\|F\| < 10^{-30}\). The last six iterations show
Newton-quadratic convergence as \(\|F\|\) collapses from \(4.2 \times 10^{-20}\)
to \(2.95 \times 10^{-34}\). From that point the \(115\)-digit Newton solve of
the full equioscillation system converges to a residual of
\(5.6 \times 10^{-98}\), and the distributed \(K = 11\) constants are its
result.

\Cref{tab:remedies} summarises the pattern. Three values of \(K\) produced
three basin geometries and three remedies, and no single global-search pipeline
succeeded in all three cases.

\begin{table}[htbp]
    \footnotesize
    \caption{Basin geometry, trap, and remedy at \(K = 9, 10, 11\). Each case
        requires a nonlocal modification of the search dynamics matched to
        its geometry. CMA-ES denotes covariance matrix adaptation evolution
        strategy, LM denotes Levenberg--Marquardt, \(F\) is the
        equioscillation residual, and \(J\) is its Jacobian.}
    \label{tab:remedies}
    \centering
    \begin{tabular}{rlll}
        \toprule
        \(K\) & Basin geometry & Trap & Remedy \\
        \midrule
        9  & wrong basin dominant & nonequioscillating stationary point
           & global escape, CMA-ES \\
        10 & right basin, algebraic saddle & \(J^{\!\top} F = 0\)
           & perturbed-reference Remez \\
        11 & right basin, narrow & full Newton step exits the basin
           & chain LM, 100 digits \\
        \bottomrule
    \end{tabular}
\end{table}

\subsection{Reproducibility}
\label{sec:reproducibility}
The certified constants for every entry in \cref{tab:cert} are publicly available
at \url{https://github.com/li-ruijie/LogisNormIntegral}, as the per-\(K\) files
\texttt{inst/certs/minimax\_k01.csv} through \texttt{minimax\_k11.csv} for the weights and
scales, and \texttt{inst/certs/deltas.csv} for the bounds. The relaxed
Bernstein--Akhiezer certificates at \(K = 14, \ldots, 17\) and the signed
\(K = 22\) frontier are in the same directory, as
\texttt{minimax\_k14\_BA\_cert.csv} through \texttt{minimax\_k17\_BA\_cert.csv}
and \texttt{minimax\_k22\_signed\_BA\_cert.csv}, with their audited bounds in
the same \texttt{deltas.csv}. Three working precisions were used in the solves. The
equioscillation audits and the \(K = 11\) chain Levenberg--Marquardt run were
carried out in 100-significant-digit arithmetic (Python
\texttt{mpmath}~\cite{johansso-mpmathp-2023}), and the Newton solves of the
\(K \le 11\) constants carried \(115\) working digits. The \(K = 9\) and \(K = 10\)
Newton polishes and the \(K = 11\) coordinate-space and log-space polishes of
\cref{sec:pipelines} used 256-bit MPFR arithmetic~\cite{fousse-mpfrmul-2007} (R
\texttt{Rmpfr}~\cite{machler-rmpfrmp-2025}). The global searches ran in double
precision. The bounds in \texttt{deltas.csv} were then recomputed by
\cref{alg:certify} at \(140\) significant digits, also in \texttt{mpmath}.
The certificates are floating-point verifications at the stated precisions rather
than interval enclosures. All cross-checks at \(50\) versus \(100\) significant
digits agree to \(50\) or more digits, placing the rounding effect far below the
reported digits.

During the preparation of this work the author used Claude Opus 5 (Anthropic) to
perform language editing and to act as a coding assistant. After using this
tool, the author reviewed and edited the content as needed and takes full
responsibility for the content of the article.

\section{Structure beyond \texorpdfstring{\(K = 11\)}{K = 11}}
\label{sec:beyond}

Above \(K = 11\) no search reached a strict \(2K\)-point equioscillation, so the
question becomes what can still be proved and certified there. This section
establishes the nonincreasing bound, states the plateau conjecture at
\(K = 12, 13\), refutes the same equality at \(K = 14, 15\), maps the
relaxed certificates through \(K = 22\), proves that \(\Delta_K^\star \to 0\)
(\cref{prop:limit}), and certifies the signed \(K = 22\) frontier.

\subsection{The nonincreasing bound and the finite plateau}

The strict \(2K\)-alternation regime of \cref{cor:equiosc} holds for
\(K = 1, \ldots, 11\). Beyond \(K = 11\) the certified minimax is nonincreasing,
\(\Delta_K^\star \le \Delta_{11}^\star\), since padding the \(11\)-atom optimum
with \(K - 11\) further atoms of total weight \(\varepsilon\) perturbs its error
by at most \(\varepsilon\), and \(\varepsilon\) is arbitrarily small.

\begin{conjecture}[Finite plateau] \label{conj:plateau}
    \(\Delta_{12}^\star = \Delta_{13}^\star = \Delta_{11}^\star\).
\end{conjecture}

Only the \(\le\) direction is proved, and the conjectured plateau cannot
extend past \(K = 13\). At \(K = 14, 15\) the corresponding equality is
refuted by exhibited mixtures audited in 100-digit arithmetic. They are
certified under a relaxed Bernstein--Akhiezer
criterion~\cite{akhiezer-theorya-1956}, written BA-\((K{+}1)\) below, which asks
for \(K + 1\) critical points, rather than \(2K\), at level at least \(0.9\)
times the audited supremum error (\cref{tab:hierarchy}). The two
witnesses place \(\Delta_{14}^\star \le 7.233 \times 10^{-12}\) and
\(\Delta_{15}^\star \le 7.211 \times 10^{-12}\), strictly below
\(\Delta_{11}^\star\). At \(K = 12, 13\) three independent \(3{,}000\)-seed
Sobol multistart streams found nothing below \(\Delta_{11}^\star\), and runs
seeded at the padded \(K = 11\) optimum returned to it. Reoptimising the
weights on a \(13\)-scale subset of each \(K = 14, 15\) witness, with the
scales held fixed, likewise gave no error below \(\Delta_{11}^\star\). These
three instruments (seeded multistart, constrained-Remez weight reoptimisation,
and scale-subset selection) are the same ones that produced the \(K = 14, 15\)
witnesses. The natural route of counting zeros is inconclusive, since the
sharp nodal bound of \cref{lem:exp-sum-zeros} permits the relevant signed
mixture at least as many interior zeros as equioscillation of the \(K = 11\)
certificate forces.

\subsection{The certified hierarchy and the cross-\texorpdfstring{\(K\)}{K} records}

\Cref{tab:hierarchy} collects what is certified from \(K = 11\) through
\(K = 17\) and at \(K = 22\), namely the strict certificate at \(K = 11\), the
relaxed certificates above it, and the two \(K = 22\) records. Mixtures with
positive weights, as in \eqref{eq:mixture}, form the positive class, and
mixtures with weights of either sign and unit sum form the signed class of
\cref{sec:frontier}.

\begin{table}[htbp]
    \footnotesize
    \caption{Certified positive-class hierarchy for \(K = 11, \ldots, 17\), the
        positive-class \(K = 22\) record (audited at higher precision, no
        certificate file shipped), and
        the best-known signed \(K = 22\) frontier, extending \cref{tab:cert}. The
        conjectured plateau at \(K = 12, 13\) is the equality
        \(\Delta_K^\star = \Delta_{11}^\star\) of \cref{conj:plateau} (only
        \(\le\) is proved). The Bernstein--Akhiezer BA-\((K{+}1)\) certificate
        relaxes the strict \(2K\)-equioscillation to \(K{+}1\) critical points,
        and the BA-\(K\) certificate to \(K\), each at level at least \(0.9\) times
        the audited supremum error of the exhibited mixture. The \(K = 11\), \(14\),
        \(15\), \(16\), \(17\), and \(22\) entries are verified in 100-digit
        arithmetic. Relaxed-certificate entries are rounded up at the last
        shown digit, so each printed value is itself a certified upper bound
        on \(\Delta_K^\star\).}
    \label{tab:hierarchy}
    \centering
    \begin{tabular}{rlll}
        \toprule
        \(K\)  & Certified error (100-digit audit) & Certificate & Source \\
        \midrule
        11     & \(1.258 \times 10^{-11}\) & strict \(2K\)-equioscillation & \cref{tab:cert} \\
        12, 13 & \(\le \Delta_{11}^\star\) & plateau conjectured (\(\le\) proved) & this work \\
        14     & \(7.233 \times 10^{-12}\) & BA-\((K{+}1)\) relaxed & this work \\
        15     & \(7.211 \times 10^{-12}\) & BA-\((K{+}1)\) relaxed & this work \\
        16     & \(4.066 \times 10^{-12}\) & BA-\((K{+}1)\) relaxed & this work \\
        17     & \(2.544 \times 10^{-12}\) & BA-\(K\) relaxed & this work \\
        \addlinespace
        22     & \(2.811 \times 10^{-14}\) & BA-\((K{+}1)\) relaxed, positive class & this work \\
        22     & \(3.605 \times 10^{-15}\) & signed BA-\((K{+}1)\) relaxed & this work \\
        \bottomrule
    \end{tabular}
\end{table}

The relaxed certificates retain the upper-bound half of the strict guarantee,
namely that the exhibited mixture's audited supremum error bounds \(\Delta_K^\star\)
from above at the stated precision, without the de la Vall\'ee Poussin lower
bound. The relaxation concerns only the equioscillation structure exhibited, and
no optimality is claimed for \(K \ge 12\). Under the BA-\((K{+}1)\) criterion,
the \(K = 14\) and \(K = 15\) mixtures refuting the plateau equality have
audited supremum errors \(7.23290 \times 10^{-12}\) and \(7.21007 \times 10^{-12}\),
and \(K = 16\) reaches \(4.065 \times 10^{-12}\). At \(K = 17\) the \(K\)-point
BA-\(K\) criterion certifies \(2.544 \times 10^{-12}\). Beyond \(K = 17\) the
Bernstein--Akhiezer criterion ceases to hold in the positive class under the
double-precision linear-programming grid pipeline, which optimises the scales
around a grid-based Chebyshev weight step (\cref{sec:pipelines}) and
whose weight step floors near
\(2.4 \times 10^{-13}\). That floor is an algorithmic artefact, not a minimax
bound. A potential-theoretic reading of the plateau, as single-interval
saturation before a bifurcation to a two-interval equilibrium measure, with
the onset of improvement at \(K = 14\), is given in \cref{sec:S-two-interval}.

\Cref{tab:critmap} records how the strict regime degrades. For
\(K = 1, \ldots, 11\) every certified minimiser places all \(2K\) critical
points at the supremum. At \(K = 12\) the count \(n_{0.99}\) of critical points
within \(1\%\) of the supremum collapses, falling to as low as one by \(K = 19\),
while the raw count \(n_{\mathrm{crit}}\) stays near order \(K\). The transition
is not a rank deficiency. A rank ceiling in the
double-precision singular value decomposition first suggested an effective
dimension near \(22\). The equioscillation system in the weights, the scales,
and the minimax level has a \(2K \times (2K+1)\) Jacobian. In \(60\)-digit
arithmetic this matrix has full row rank at \(K = 16\) (\(32/32\)), \(K = 17\)
(\(34/34\)), and the \(K = 22\) positive-class record (\(44/44\), smallest
singular value \(\approx 1.3 \times 10^{-21}\)). The strict \(2K\) system is
therefore numerically well posed, and the obstruction is finding the global
basin.

\begin{table}[htbp]
    \footnotesize
    \caption{Cross-\(K\) critical-point map of the positive-class pipeline
        record. The column \(n_{\mathrm{crit}}\) counts distinct critical points
        of the error of the tabulated configuration on \(0 < z \le 200\), and
        \(n_{0.99}\) those
        within \(1\%\) of the supremum, and strict equioscillation requires
        \(n_{0.99} = 2K\). Past that window the logistic's exponential tail
        dominates every Gaussian tail, so the error is monotone there except at
        \(K = 19, 20, 21\), each of which carries one further extremum, near
        \(z = 486\), \(294\), and \(387\) at magnitudes below \(10^{-128}\). The
        \(K \ge 12\) rows are the record of the double-precision
        linear-programming grid pipeline. Its direct nondegenerate solutions at
        \(K = 12, \ldots, 15\) have supremum errors above the padding
        bound \(\Delta_K^\star \le \Delta_{11}^\star\) of the text and above
        the certified \(K = 14, 15\) bounds of \cref{tab:hierarchy}, so the map
        records critical-point structure rather than best bounds. At \(K = 16\)
        the supremum error is the audited certificate value, while the two count
        columns give ranges spanning the near-best pipeline configurations. The
        \(K = 22\) row is superseded in the positive class by the higher-precision
        entry of \cref{tab:hierarchy} (\(2.811 \times 10^{-14}\)) and in the
        signed class by the frontier of the text (30 dominant alternations at
        \(3.604 \times 10^{-15}\)).}
    \label{tab:critmap}
    \centering
    \begin{tabular}{rrrl}
        \toprule
        \(K\) & \(n_{\mathrm{crit}}\) & \(n_{0.99}\) & Supremum error \\
        \midrule
        1--11 & \(2K\) & \(2K\) & \cref{tab:cert}, strict \\
        12 & 16 & 15 & \(5.04 \times 10^{-10}\) \\
        13 & 20 & 15 & \(1.57 \times 10^{-10}\) \\
        14 & 18 & 15 & \(6.30 \times 10^{-11}\) \\
        15 & 19 & 13 & \(1.72 \times 10^{-11}\) \\
        16 & 21--25 & 4--16 & \(4.07 \times 10^{-12}\) \\
        17 & 21 & 3 & \(2.54 \times 10^{-12}\) \\
        18 & 22 & 2 & \(4.23 \times 10^{-13}\) \\
        19 & 25 & 1 & \(3.30 \times 10^{-13}\) \\
        20 & 21 & 2 & \(2.97 \times 10^{-13}\) \\
        21 & 20 & 1 & \(2.52 \times 10^{-13}\) \\
        22 & 23 & 1 & \(2.41 \times 10^{-13}\) \\
        \bottomrule
    \end{tabular}
\end{table}

\subsection{The vanishing of the minimax error in the limit}
\label{sec:limit}

The padding bound above caps \(\Delta_K^\star\) by \(\Delta_{11}^\star\) at
every \(K \ge 12\), leaving open whether the sequence descends to zero at all.
An explicit construction that discretises the mixing density settles it.

\begin{proposition}[Vanishing minimax error] \label{prop:limit}
    \(\lim_{K \to \infty} \Delta_K^\star = 0\).
\end{proposition}

\begin{proof}
    Define the scale-invariant modulus
    \(D(\rho) = \sup_{u > 0} |\Phi(\rho u) - \Phi(u)|\). For \(\rho \ge 1\),
    monotonicity of \(\phi\) on \((0, \infty)\) gives
    \(|\Phi(\rho u) - \Phi(u)| = \int_u^{\rho u} \phi \le (\rho - 1)\, u \phi(u)
    \le (\rho - 1)/\sqrt{2 \pi e}\), since \(\max_{u > 0} u \phi(u) = \phi(1)\),
    and the case \(\rho \le 1\) is symmetric, so
    \(\widetilde D(\rho) = \sup_{r \in [1/\rho, \rho]} D(r) \to 0\) as
    \(\rho \downarrow 1\).

    Fix \(K\) and partition \([1/K, K]\) logarithmically
    into \(K\) intervals \(I_j\) of common ratio \(\rho_K = K^{2/K} \to 1\). Place
    an atom \(s_j \in I_j\) of weight \(w_j = \int_{I_j} q\, ds\), adding the tail
    mass \(\tau_K = \int_{(0, 1/K) \cup (K, \infty)} q\, ds\) to the end atoms, so
    \(\sum_j w_j = 1\) and the mixture lies in \(\F_K\). For \(s \in I_j\) and
    \(z > 0\), writing \(u = s_j z\) and \(\rho = s / s_j\) gives
    \(|\Phi(s z) - \Phi(s_j z)| \le \widetilde D(\rho_K)\). By the scale-mixture
    representation \(\Lam(z) = \int_0^\infty \Phi(s z)\, q(s)\, ds\) of
    \cref{sec:intro},
    \[
        \Bigl| \Lam(z) - \sum_{j=1}^K w_j \Phi(s_j z) \Bigr|
        \le \sum_j \int_{I_j} |\Phi(s z) - \Phi(s_j z)|\, q\, ds + 2 \tau_K
        \le \widetilde D(\rho_K) + 2 \tau_K,
    \]
    where \(2 \tau_K\) bounds both tail contributions. The error is odd in \(z\),
    by \(\Lam(-z) = 1 - \Lam(z)\) and \(\sum_j w_j = 1\), so the bound holds on
    all of \(\mathbb{R}\). As \(K \to \infty\), \(\widetilde D(\rho_K) \to 0\) and
    \(\tau_K \to 0\), so \(\Delta_K^\star \to 0\).
\end{proof}

The natural limiting object is the mixing density \(q\) itself, available in
closed form as \(q(s) = (2/s^3) \sum_{n \ge 1} (-1)^{n-1} n^2
\exp(-n^2 / (2 s^2))\)~\cite{stefansk-normals-1991}. Its low-order moments admit
closed forms through the Dirichlet eta function, catalogued in
\cref{sec:S-moments}, and the certified atoms already realise two of them. The
second moment \(\sum_i w_i^\star (s_i^\star)^2\) approaches
\(7 \zeta(3) / (2 \pi^2) = 0.42627\ldots\) (seven digits by \(K = 10\)), and the
log-scale moment \(\sum_i w_i^\star \log s_i^\star\) approaches
\((\gamma_{\mathrm{EM}} + \log 2 - 2 \log \pi)/2 = -0.50954\ldots\) (seven
digits), with \(\gamma_{\mathrm{EM}}\) the Euler--Mascheroni constant.

Both limits first appeared as byproduct convergents of a search for a closed
form of the certified constant \(\Delta_{11}^\star\) itself. Integer-relation
searches (PSLQ)~\cite{ferguson-analysis-1999} returned no relation in three
campaigns. The first ran over the modular-transform constants of the
associated Mordell integral~\cite{zwegers-mockthe-2002} at coefficient heights
up to \(10^{50}\). The second ran over \(81\) statistics of the certified
weights and scales in nine families, together with \(21\) power-law fits of the
rank-ordered scales, at heights up to \(10^{25}\). The third ran over the
limiting values of the statistics that converge with \(K\), the two limits
above among them, at heights up to \(10^{60}\).
The certified constant therefore satisfies no integer relation of low height
at the precision tested, and the two convergents are moments of the limiting
mixing density \(q\) rather than constants special to \(K = 11\).

\subsection{The signed frontier at higher precision}
\label{sec:frontier}

Working at higher precision removes the double-precision floor. A scale search
in the positive class whose objective rewards the alternation count reaches
\(2.810 \times 10^{-14}\) at \(K = 22\) under the BA-\((K{+}1)\) criterion, the
corresponding entry of \cref{tab:hierarchy}. A covariance matrix adaptation
evolution strategy (CMA-ES) run with the same objective in the signed class
(\(\sum_i w_i = 1\), all \(w_i \ne 0\), signs unrestricted) reaches
\(\Delta = 3.604 \times 10^{-15}\) at \(K = 22\), whose
error has \(30\) strictly alternating extrema at level \(\ge 0.99\,\|e\|_\infty\).
The bound lies within a factor of about thirty-three of the double-precision unit
roundoff \(2^{-53}\).

Solving the prescribed-alternation Newton system in normal-equation form stalls
near \(4.5 \times 10^{-20}\), an artefact of squaring the Jacobian condition
number (\(\kappa \approx 1.8 \times 10^{15}\) becomes \(\kappa^2 \approx
3 \times 10^{30}\), exhausting the working precision of that solver stage).
Stepping instead on the singular value decomposition (SVD) of \(J\) directly,
with per-direction damping and
geodesic acceleration~\cite{transtru-whyare-2010, transtru-geometry-2011},
drives the residual to \(1.05 \times 10^{-90}\) and certifies a strict signed
\(30\)-point equioscillation. The error at this frontier has \(31\) critical
points, the \(30\) dominant alternations and one more below the level.

A certificate exact at double precision, one placing the error below the
double-precision unit roundoff \(2^{-53}\), would require the full strict
\(44\)-point equioscillation. That count is unreachable from this frontier by
any route attempted, and the diagnostics agree on the mechanism. A function
homotopy from an explicit \(44\)-lobe equiripple seed is conditioning-blocked,
since the seed carries \(\kappa \approx 2 \times 10^{97}\). Growing the atom count
to \(K = 23\) and \(K = 24\) either returns to the frontier or reaches
basins six orders of magnitude worse, with cancelling weights. Seeding that
SVD-based solver within the \(K = 22\) basin with the missing critical points
inserted stalls at \(\|F\| \approx 3.6 \times 10^{-14}\) to
\(5.9 \times 10^{-14}\) and returns to the frontier's \(31\) critical points.
The seed conditioning (\(\kappa \approx 10^{28}\) to \(10^{35}\)) places the
attainable floor near \(10^{-45}\) to \(10^{-52}\) at that solver's working
precision,
so the stall is structural rather than a precision limit. Creating the \(13\) missing
critical points is a bifurcation through
degenerate configurations with \(e' = e'' = 0\), where the corresponding
Jacobian rows become linearly dependent, and a local Newton flow does not cross such a fold.

A signed mixture bounds the integral by the same triangle inequality as a
positive one, \(|I(\eta, v) - \hat I(\eta, v)| \le \Delta\) with \(\hat I\) as in
\eqref{eq:int-from-mixture}, whatever the signs of the weights. The mediation
application does not depend on these higher-\(K\) facts, since the certified
\(K = 11\) bound of \(1.26 \times 10^{-11}\) is already far below the sampling
variation of the fitted coefficients.

\section{Application to causal mediation and comparison with alternative
    evaluators}
\label{sec:mediation}

\Cref{sec:certification,sec:beyond} certified the mixture constants and stated
\cref{alg:evaluate} for evaluating the logistic-normal integral
\eqref{eq:logistic-normal} to a fixed uniform error bound. This section applies
that evaluation to causal mediation with a binary outcome, where the natural
indirect effect reduces to a difference of two logistic-normal integrals.
The section compares the certified mixture evaluation with the rare-outcome
closed form of Valeri and VanderWeele, across the strength of the mediator and
on public data (\cref{sec:rare-outcome}), and forms quasi-Bayesian
intervals (\cref{sec:qb-comparison}). It then benchmarks
the mixture against Gauss--Hermite quadrature, Gauss--Fermi approximations,
rational pole expansions, and Pirjol's Poisson-summation series
(\cref{sec:pointwise}), and shows that the closed form composes across \(k\)
sequential or parallel normal mediators at the cost of a single evaluation
(\cref{sec:multi}).

\subsection{The integral in causal mediation}
Consider a two-stage model for a binary outcome \(Y \in \{0, 1\}\) and a
continuous mediator \(M\) given exposure \(z\) and covariates \(c\). The
outcome follows the logistic regression
\(P(Y = 1 \mid z, m, c) = \Lam(\theta_0 + \theta_1 z + \theta_2 m + \theta_3 z m +
\theta_4^{\!\top} c)\), and the mediator follows the homoscedastic normal
regression \(M \mid z, c \sim N(\gamma_0 + \gamma_1 z + \gamma_2^{\!\top} c,
\sigma_M^2)\). Under the standard identification
assumptions and cross-world independence~\cite{vanderwe-explanat-2015}, the
expected outcome when the direct exposure is held at \(z\) and the mediator is
drawn under exposure \(z^\prime\) is the expectation of
\(\Lam(\theta_0 + \theta_1 z + (\theta_2 + \theta_3 z) M + \theta_4^{\!\top} c)\)
over \(M \sim N(\gamma_0 + \gamma_1 z^\prime + \gamma_2^{\!\top} c,
\sigma_M^2)\). That expectation is the logistic-normal integral
\eqref{eq:logistic-normal} with
\(\eta(z, z^\prime, c) = \theta_0 + \theta_1 z + (\theta_2 + \theta_3 z)(\gamma_0 +
\gamma_1 z^\prime + \gamma_2^{\!\top} c) + \theta_4^{\!\top} c\) and
\(v(z) = (\theta_2 + \theta_3 z)^2 \sigma_M^2\). The natural indirect
effect of an exposure shift \(0 \to z^\prime\) in the mediator is
the counterfactual contrast
\(\text{NIE}(z, z^\prime, c) = \E[Y(z, M(z^\prime)) - Y(z, M(0)) \mid c] =
I(\eta(z, z^\prime, c), v(z)) - I(\eta(z, 0, c), v(z))\), and setting
\(z = z^\prime = 1\) gives the total natural indirect effect. Normality of \(M\)
reduces each potential-outcome mean to a single integral
\eqref{eq:logistic-normal}, and homoscedasticity makes \(v\) independent of
\(z^\prime\) and \(c\).

The closed forms of Valeri and VanderWeele~\cite{valeri-mediatio-2013}, the
default in applied causal mediation through their SAS and SPSS macros, avoid
the integral by assuming that the outcome is rare. The logistic function is
then close to the exponential, the expectation over the normal mediator is a
moment generating function~\cite[Proposition~2.4]{vanderwe-explanat-2015}, and
the effects follow in closed form on the odds-ratio scale. For the total
natural indirect effect, \(z = z^\prime = 1\), the form is
\begin{equation} \label{eq:rare-form}
    \mathrm{OR}^{\text{NIE}}_{\text{rare}} = \exp\{(\theta_2 + \theta_3)\, \gamma_1\},
\end{equation}
which involves neither the intercept \(\theta_0\), nor the covariates, nor
\(\sigma_M\). The quantity it approximates is the odds ratio of the two
potential-outcome means,
\begin{equation} \label{eq:exact-or}
    \mathrm{OR}^{\text{NIE}}(c) =
    \frac{\operatorname{odds} I(\eta(1, 1, c), v(1))}{\operatorname{odds} I(\eta(1, 0, c), v(1))},
    \qquad \operatorname{odds} p = \frac{p}{1 - p},
\end{equation}
and \eqref{eq:rare-form} is the limit of \eqref{eq:exact-or} as
\(\theta_0 \to -\infty\), where \(\operatorname{odds} I(\eta, v) \sim I(\eta, v)
\sim e^{\eta + v/2}\) and the factors common to \(\eta(1, 1, c)\) and
\(\eta(1, 0, c)\) cancel. The
forms are accurate when the outcome is rare and biased when it is
common~\cite{gaynor-mediatio-2019}. The certified mixture
\eqref{eq:int-from-mixture} evaluates both integrals of \eqref{eq:exact-or} to
within \(\Delta_K^\star\) at every \((\eta, v)\), so it gives the odds ratio,
and the effect on the probability scale, without the assumption.
\Cref{sec:rare-outcome} measures how far the two answers lie apart.

Deployed software either evaluates the same integral or avoids it, by one of
three methods.
Quasi-Bayesian simulation, as in the
\texttt{mediation} package~\cite{tingley-mediatio-2014}, draws the mediator and
averages the fitted outcome probability, carrying Monte Carlo error and the
cost of an inner sampling loop. The rare-outcome closed forms, implemented in
\texttt{regmedint}~\cite{yoshida-regmedin-2026}, avoid the integral and hold
when the outcome is rare. Fixed-grid numerical integration,
the route taken by Mplus~\cite[Section~8.1.3]{muthen-regressi-2016}, evaluates the integrand on
a fixed \(101\)-point grid over the standardised range \([-5, 5]\), a
deterministic rule whose accuracy is fixed by the grid spacing.
\Cref{sec:S-comparison} compares the mixture with simulation in accuracy and
runtime, and relative to the \(101\)-node grid the mixture evaluates a shorter
sum with a certified error bound.

\subsection{The rare-outcome closed form and a worked example}
\label{sec:rare-outcome}

Consider a synthetic instance with \((\theta_0, \ldots, \theta_4) = (-1.0, 0.5,
1.5, 0.5, -0.2)\), \((\gamma_0, \gamma_1, \gamma_2) = (0.5, 1.2, 0.3)\), and
\(\sigma_M = 0.8\), evaluated at \(z = z^\prime = 1\) and \(c = 0\) so that
\(v = 2.56\). On this instance, the certified \(K = 11\) mixture yields
\(\text{NIE} = 0.30464\) with absolute error at most \(2\Delta_{11}^\star <
2.52 \times 10^{-11}\), twice the per-integral bound since the effect is a
difference of two integrals. The two potential-outcome means are \(0.8896\) and
\(0.5849\), so the exact odds ratio \eqref{eq:exact-or} is \(5.716\). The closed
form \eqref{eq:rare-form} gives \(e^{2.4} = 11.02\), larger by a factor of
\(1.93\).

That factor is specific to this instance, which has a common outcome and a
strong mediator, since
\(v = 2.56\) means that one residual standard deviation of the mediator
multiplies the odds of the outcome by \(e^{1.6} \approx 5\).
\Cref{tab:rare-grid} therefore varies both. It scales the two mediator
coefficients \((\theta_2, \theta_3)\) of the instance by a common factor between
\(1/8\) and \(1\), which moves \(v\) from \(0.04\) to \(2.56\) and the closed
form from \(1.35\) to \(11.02\). For each \(v\), it sets \(\theta_0\) so that
the outcome probability \(I(\eta(1, 0, 0), v)\), when the mediator is drawn
under no exposure, takes the column value. At \(v \le 0.16\) the closed form
lies within \(2.2\%\) of the exact odds ratio in every column, including an
outcome probability of one half. At
\(v = 0.64\) it overstates the odds ratio by \(1.8\%\) when that probability is
\(0.01\) and by \(15\%\) when it is \(0.25\) or more. At \(v = 2.56\) it
overstates by \(33\%\) even when the probability is \(0.01\), since drawing
the mediator under exposure then raises the outcome probability to \(0.077\).

\begin{table}[htbp]
    \footnotesize
    \caption{The rare-outcome closed form \eqref{eq:rare-form} divided by the
        exact odds ratio \eqref{eq:exact-or} of the natural indirect effect, by
        the variance argument \(v\) and by the outcome probability when the
        mediator is drawn under no exposure. The coefficients are those of the synthetic instance
        of \cref{sec:rare-outcome} with \((\theta_2, \theta_3)\) scaled by
        \(1/8\), \(1/4\), \(1/2\), \(3/4\), and \(1\), and \(\theta_0\) solved
        for the stated probability. Reproducible via
        \texttt{inst/reproduction/Example1.R}.}
    \label{tab:rare-grid}
    \centering
    \begin{tabular}{rrrrrrr}
        \toprule
        & & \multicolumn{5}{c}{Outcome probability \(I(\eta(1, 0, 0), v)\)} \\
        \cmidrule(lr){3-7}
        \(v\) & Closed form & \(0.01\) & \(0.05\) & \(0.10\) & \(0.25\) & \(0.50\) \\
        \midrule
        \(0.04\) & \(1.350\)  & \(1.000\) & \(1.001\) & \(1.001\) & \(1.002\) & \(1.003\) \\
        \(0.16\) & \(1.822\)  & \(1.001\) & \(1.006\) & \(1.011\) & \(1.019\) & \(1.022\) \\
        \(0.64\) & \(3.320\)  & \(1.018\) & \(1.067\) & \(1.103\) & \(1.150\) & \(1.151\) \\
        \(1.44\) & \(6.050\)  & \(1.100\) & \(1.269\) & \(1.367\) & \(1.473\) & \(1.451\) \\
        \(2.56\) & \(11.023\) & \(1.330\) & \(1.709\) & \(1.902\) & \(2.089\) & \(2.010\) \\
        \bottomrule
    \end{tabular}
\end{table}

A model fitted to public data lies near the first row of that grid. In the
Western Collaborative Group Study (WCGS)~\cite{rosenman-coronary-1975}, the exposure is a
body-mass index of \(25\) or more, the mediator is systolic blood pressure, and
the outcome is coronary heart disease, with adjustment for age, cigarettes
smoked, and behaviour pattern. Of the \(3154\) participants \(8.1\%\) have the
event, and the fitted model has \(v(1) = 0.032\). The closed form gives
\(1.0826\), and the exact odds ratio \eqref{eq:exact-or}, which varies with the
covariates between \(1.0819\) and \(1.0825\), has geometric mean \(1.0824\) over
the participants. Body-mass index and blood pressure are states rather than
interventions, and both were measured at the same visit, so this computation
evaluates the mediation formula at the fitted models and offers no causal
reading of the data. The data are those of the R package \texttt{faraway}, and
over every potential-outcome mean the mixture agrees with adaptive quadrature to
\(1.6 \times 10^{-15}\) (\texttt{inst/reproduction/PublicData.R}).

The closed form is therefore accurate where \(v\) is small, as in WCGS, and
where the outcome is rare under both laws of the mediator. Where neither holds
it departs from the exact value, by \(15\%\) at \(v = 0.64\) and an outcome
probability of \(0.25\) in \cref{tab:rare-grid}, and by a factor of \(1.93\) in
the synthetic instance. It cannot report which case applies, since it uses
neither \(\theta_0\) nor \(\sigma_M\). The certified mixture returns the exact
odds ratio in each case at the cost of \(2K\) values of \(\Phi\).

The default mixture in the accompanying repository is the certified \(K = 22\) signed frontier, with
worst-case error \(3.605 \times 10^{-15}\). The positive-class \(K = 11\) mixture is displayed here and in
\cref{tab:comparison}, and the two agree to ten decimals on identical parameter draws, so the choice of
certified mixture is immaterial to the reported values.

\subsection{Quasi-Bayesian intervals}
\label{sec:qb-comparison}

Since the maximum-likelihood estimators are asymptotically normal, drawing
\(R = 10^4\) parameter vectors from their fitted sampling distribution and evaluating
\eqref{eq:int-from-mixture} at each draw gives asymmetric quasi-Bayesian
intervals~\cite{king-makingm-2000} (asymptotically valid credible limits under
the Bernstein--von Mises correspondence). Each draw evaluates two integrals at
\(K\) calls to \(\Phi\) apiece, so \(R = 10^4\) draws require \(2 K R\) values
of \(\Phi\) and complete in under \(15\) milliseconds.

\Cref{sec:S-comparison} compares three inner-integral strategies across three
outcome-prevalence regimes at \(R = 10^4\) outer draws, with the coefficients
of the synthetic instance and \(\theta_0 \in \{-3, -1, 0\}\). The certified
mixture, inner Monte Carlo, and the pointwise evaluator of \cref{sec:pointwise}
agree in every regime. Timed at the repository default \(K = 22\), the
mixture is more than \(380\) times faster than inner Monte Carlo at
\(M_{\text{inner}} = 10^4\) and more than \(13\) times faster than the pointwise
route. In a repeated-sampling experiment of \(2000\) data sets per regime,
the mixture intervals for the natural indirect effect cover its population
value in \(95.8\%\), \(95.7\%\), and \(95.2\%\) of the replicates. For the
exact odds ratio \eqref{eq:exact-or}, the mixture intervals cover the
population value in \(95.8\%\), \(95.8\%\), and \(95.2\%\), whereas intervals
from the closed form \eqref{eq:rare-form} on the same draws cover it in
\(39.6\%\), \(70.4\%\), and \(87.8\%\). These regimes share
\(v = 2.56\) and a common outcome, the last row of \cref{tab:rare-grid}, where
the closed form is furthest from the quantity it approximates.

The certified mixture therefore replaces inner Monte Carlo and fixed-grid
quadrature with a single closed-form evaluation carrying a certified uniform
error bound, at the cost of a delta-method interval, and it requires no
assumption on the outcome probability.

\subsection{Alternative evaluations and the pointwise route}
\label{sec:pointwise}
Beyond the simulation, rare-outcome, and grid methods deployed in mediation software, the
numerical and physics literatures evaluate \eqref{eq:logistic-normal} by
quadrature rules, Gauss--Fermi approximations, rational pole expansions, and
Pirjol's pointwise recursion and Poisson-summation
series~\cite{pirjol-logistic-2013} (\cref{tab:alternatives}). Among
quadrature rules, Crouch and Spiegelman~\cite{crouch-evaluati-1990} designed
one for integrals of this form with logistic-normal models as the application,
a trapezoidal rule whose step is chosen from an explicit error bound.
Standard Gauss--Hermite quadrature integrates against the normal weight
directly, but in the quadrature variable the integrand's poles at
\(\eta + m = \pm i\pi\) lie at distance \(\pi/\sqrt{2v}\) from the real axis,
so convergence degrades as \(v\) grows. With \(11\) nodes, matching the
\(K = 11\) mixture in function evaluations, its largest error over \(801\)
equally spaced values of \(\eta \in [-10, 10]\) is \(1.8 \times 10^{-7}\) at
\(v = 1\) and \(2.2 \times 10^{-3}\) at \(v = 9\). Even \(100\) nodes leave an
error of \(2.2 \times 10^{-9}\) at \(v = 9\), compared with the mixture's
certified bound \(1.26 \times 10^{-11}\) at every \(v\)
(\texttt{inst/reproduction/Quadrature.R}).

The semiconductor literature cited in \cref{sec:intro} treats the same
integral under its Gauss--Fermi name. Its two-branch
approximation~\cite{paasch-chargec-2010} is a closed form whose error the
authors put below \(10\%\) to \(15\%\) for \(-v < \eta < 0\). Measured
against adaptive quadrature over \(|\eta| \le 30\) and
\(v \in \{4, 9, 16, 25, 36\}\), widths within the range the article plots, its
largest absolute error is \(4.5 \times 10^{-3}\) and its largest relative error
is \(18\%\). The error-function
series~\cite{mehmetog-analytic-2014, selvaggi-analytic-2018} is exact in the
limit and converges slowly. It needed between \(56\) and \(492\) terms for an
error of \(10^{-6}\) at the six arguments \(\eta \in \{1, 3\}\),
\(v \in \{1, 9, 25\}\). Selvaggi presents the series as a check on numerical
methods and notes that numerical integration is the more efficient route to
values. An analytical approximation also exists for the inverse relation, the
Fermi level as a function of the carrier
density~\cite{sharma-analytic-2013}.

The moment expansion of \cref{tab:alternatives} replaces \(I(\eta, v)\) by a
truncated Taylor series of \(\Lam(\eta + m)\) around \(m = 0\), whose second-order
form is \(\Lam(\eta) + \tfrac{1}{2} v \Lam''(\eta)\). Since
\(\Lam(\eta + m)\) has poles in the complex plane at \(\eta + m = \pm i\pi\),
its Taylor series in \(m\) has the finite convergence radius
\(\sqrt{\eta^2 + \pi^2}\), so term-by-term expectation over the Gaussian
density on \((-\infty, \infty)\) diverges for every \(v > 0\) when
\(\eta \ne 0\). For the second-order truncation, the odd Gaussian moments
vanish and the fourth moment \(3v^2\) supplies the factor \(3/4! = \tfrac18\),
bounding the remainder by
\begin{equation} \label{eq:taylor-remainder}
    |R_{\text{T}}(\eta, v)| \le \tfrac{1}{8}\, v^2\, \sup_m |\Lam^{(4)}(\eta + m)|
    \le \tfrac{v^2}{8}\, C_4, \qquad C_4 = \sup_x |\Lam^{(4)}(x)| \approx 0.128.
\end{equation}
The bound is finite for every \((\eta, v)\) but grows quadratically in \(v\),
from \(0.016\) at \(v = 1\) to \(0.26\) at \(v = 4\).

A further group approximates the sigmoid before integrating. In
electronic-structure calculations the reflected logistic function
\(\Lam(-x) = 1/(1 + e^{x})\) is the Fermi--Dirac distribution, approximated by
rational functions with simple poles. Lin, Lu, Ying, and
E~\cite{lin-pole-bas-2009} construct pole expansions whose cost grows
logarithmically with the width of the approximation interval.
Moussa~\cite{moussa-minimax-2016} tabulates the minimax rational approximant on
a half-line in double precision and notes that the conventional Remez algorithm
cannot compute it when the number of poles is large. Since a finite sum of
simple poles vanishes at infinity whereas \(\Lam(-x)\) tends to one as
\(x \to -\infty\), the error of that family tends to one there, and no member
approximates \(\Lam(-x)\) uniformly on the whole real line. A normal-scale
mixture shares both asymptotic limits of \(\Lam\), which makes a uniform bound
on \(\R\) possible. Within the closed-form rows of
\cref{tab:alternatives}, probit scaling \(\Phi(c \eta / \sqrt{1 + c^2 v})\) is
the \(K = 1\) member of \eqref{eq:int-from-mixture}. Gaynor, Schwartz, and
Lin~\cite{gaynor-mediatio-2019} used it for closed-form mediation point
estimates under a common outcome, at the two-digit accuracy floor that
\(\Delta_1^\star\) makes precise. Gory, Craigmile, and
MacEachern~\cite{gory-classge-2021} evaluate the integral by the eight-term
mixture of Monahan and Stefanski and apply Pirjol's recursion first, so that
the mixture is evaluated near the origin. Among those closed-form rows, the
normal-scale mixture family alone combines a closed form in \((\eta, v)\) with
a parameter-free certified bound.

\begin{table}[htbp]
    \footnotesize
    \caption{Alternative evaluations of \eqref{eq:logistic-normal} beyond those
        in \cref{sec:qb-comparison} and \cref{sec:S-comparison}. The probit row
        is the \(K = 1\) member of the mixture family \eqref{eq:mixture}, whose
        optimised constant attains \(\Delta_1^\star\) of \cref{tab:cert}. The
        Laplace contraction threshold follows from \(\sup_x \Lam'(x) = 1/4\), and
        the Pad\'e range keeps the sigmoid argument within the accuracy region
        \(|x| \lesssim 3\) for two standard deviations of the normal argument.
        The Gauss--Hermite figure is the largest error over \(801\) equally
        spaced values of \(\eta \in [-10, 10]\). The figures of the two-branch approximation and of
        the error-function series are measured on the arguments given in the
        text. The pointwise route is treated in the text.}
    \label{tab:alternatives}
    \centering
    \setlength{\tabcolsep}{4pt}
    \begin{tabular}{@{}L{3.5cm} L{1.8cm} L{3.5cm} L{3.4cm}@{}}
        \toprule
        Route & Closed form & Uniform guarantee & Limitation \\
        \midrule
        Probit scaling \(\Phi(c \eta / \sqrt{1 + c^2 v})\) & yes
            & \(2.3 \times 10^{-2}\) at \(c = \sqrt{3}/\pi\), \(\Delta_1^\star\) optimised
            & two-digit accuracy floor \\
        \addlinespace[2pt]
        Taylor (moment) expansion & yes
            & none, remainder \eqref{eq:taylor-remainder}
            & divergent for every \(v > 0\) when \(\eta \ne 0\), radius
              \(\sqrt{\eta^2 + \pi^2}\) \\
        \addlinespace[2pt]
        Laplace & implicit mode solve
            & none
            & skew error at moderate \(v\), contraction needs \(v < 4\) \\
        \addlinespace[2pt]
        Pad\'e in \(\tanh\) & yes
            & none
            & wrong tails, needs \(|\eta| + 2\sqrt{v} \lesssim 3\) \\
        \addlinespace[2pt]
        Variational minorants~\cite{jaakkola-bayesian-2000, anceschi-optimal-2026} & yes
            & one-sided
            & lower bounds rather than evaluations \\
        \addlinespace[2pt]
        Gauss--Hermite quadrature, \(n\) nodes & yes
            & none, error grows with \(v\)
            & \(2.2 \times 10^{-3}\) at \(n = 11\), \(v = 9\) \\
        \addlinespace[2pt]
        Two-branch approximation~\cite{paasch-chargec-2010} & yes
            & none
            & \(4.5 \times 10^{-3}\) absolute, \(18\%\) relative \\
        \addlinespace[2pt]
        Error-function series~\cite{mehmetog-analytic-2014, selvaggi-analytic-2018}
            & infinite series
            & none, no truncation bound
            & \(56\) to \(492\) terms for \(10^{-6}\) \\
        \addlinespace[2pt]
        Pirjol recursion and Poisson series~\cite{pirjol-logistic-2013} & per-point routine
            & per-point truncation control
            & no closed form in \((\eta, v)\) \\
        \addlinespace[2pt]
        Certified mixture \eqref{eq:int-from-mixture} & yes
            & \(\Delta_K^\star\), \(1.26 \times 10^{-11}\) at \(K = 11\)
            & fixed approximation floor per \(K\) \\
        \bottomrule
    \end{tabular}
\end{table}

The pointwise route of Pirjol~\cite{pirjol-logistic-2013} stands apart from
every closed-form row. The exact recursion and the Poisson-summation series
compute a single value of \eqref{eq:logistic-normal} to the accuracy of the
working arithmetic.

The comparison rests on measurement, and the script \texttt{PirjolTiming.R}
is distributed in \texttt{inst/reproduction} with its captured output alongside. On an
eight-point grid spanning \(v \le 40\), its R implementation at
truncation order \(30\) reproduces Gauss--Kronrod reference values to at most
\(2.3 \times 10^{-16}\) at \(21.5\,\mu\mathrm{s}\) per call, and at order
\(10\) to at most \(1.6 \times 10^{-12}\) at \(8.7\,\mu\mathrm{s}\). The
certified \(K = 11\) sum evaluates the same points to at most
\(2.4 \times 10^{-14}\) observed error, within its certified bound, at
\(1.7\,\mu\mathrm{s}\) per call. Per-point accuracy therefore favours the
pointwise route at sufficient truncation order, and both routes are
microsecond-scale.

What the pointwise route does not supply is the object the mediation
application consumes. Each call returns a number rather than a function, so the closed-form effect
expressions, their derivatives, and array evaluation across the \(10^4\)
quasi-Bayesian draws of \cref{sec:qb-comparison} are unavailable. In addition,
the truncation order must be controlled at every draw, since the series'
convergence rate varies with \((\eta, v)\). A subsequent
treatment~\cite{pirjol-logistic-2026} singles out an optimal series among the
possible expansions and supplies tail error control, which sharpens that
per-draw management without changing the per-point character of the guarantee.
A contemporaneous
treatment of the related logit-normal moments~\cite{holmes-inferent-2026}
illustrates the same trade-off from the other side. It approximates the logistic
function by a series truncated at eleven terms, with a fitted near-zero
correction, overflow guards on the exponential terms, and a regime threshold
selecting between branches. It reports accuracy as a measured range rather
than a bound. A certified mixture
carries the a priori uniform guarantee \(\Delta_K^\star\) at every point,
fixed before any draw is taken. The two routes are therefore complementary.
For individual values at known parameters, Pirjol's evaluators are
preferable. For the integral consumed as a function inside estimands,
sampling loops, and proofs, the certified mixture is preferable.

Inside the quasi-Bayesian loop the pointwise route is a usable inner
engine. It reproduces the certified mixture to all six decimals the script
carries in the low- and moderate-prevalence regimes, at \(0.17\) seconds per
\(10^4\) draws at truncation order \(10\) and \(0.43\) at order \(30\). It
nonetheless needs the truncation order managed at every draw and numerical
differentiation for derivatives, and it returns
nonfinite values at both edges of its domain, which the high-prevalence regime
reaches. \Cref{sec:S-comparison} measures each of these costs. Statistical noise
dominates the numerical error of both routes by at least eight orders of
magnitude, and the mixture's advantage in this application is therefore the
robustness of the guarantee and the closed form rather than the digits.

\subsection{Multiple mediators: Composition at unchanged cost}
\label{sec:multi}
The closed form composes. Identification of nested multiple-mediator
counterfactuals requires assumptions beyond the single-mediator set of
\cref{sec:mediation}~\cite{vanderwe-explanat-2015}, and this subsection's
claims are computational, valid wherever the nested mean is identified as the
stated Gaussian functional. Take two sequential mediators with outcome model
\(\operatorname{logit} P(Y = 1 \mid z, m_1, m_2) = \theta_0 + \theta_1 z +
\theta_2 m_1 + \theta_3 m_2\), where exposure-mediator interactions are permitted
since \(z\) is held fixed inside each potential outcome, and normal mediator
models \(M_1 \mid z'' \sim N(\gamma_{10} + \gamma_{11} z'', \sigma_1^2)\) and
\(M_2 \mid z', m_1 \sim N(\gamma_{20} + \gamma_{21} z' + \gamma_{22} m_1,
\sigma_2^2)\). The nested two-dimensional potential-outcome mean collapses
exactly,
\begin{equation} \label{eq:collapse}
    \E\bigl[Y\bigl(z,\, M_1(z''),\, M_2(z', M_1(z''))\bigr)\bigr]
        = I(\eta_{\mathrm{tot}},\, v_{\mathrm{tot}}),
\end{equation}
with \(c_1 = \theta_2 + \theta_3 \gamma_{22}\) the total coefficient of
\(m_1\) on the logit scale, \(\eta_{\mathrm{tot}} = \theta_0 + \theta_1 z +
\theta_3 (\gamma_{20} + \gamma_{21} z') + c_1 (\gamma_{10} + \gamma_{11}
z'')\), and \(v_{\mathrm{tot}} = \theta_3^2 \sigma_2^2 + c_1^2 \sigma_1^2\).
The argument is the heat semigroup, the reading of
\eqref{eq:logistic-normal} that Pirjol~\cite{pirjol-logistic-2013} develops. As
a function of \((\eta, v)\), the integral \(I(\eta, v) = (\Lam * \phi_v)(\eta)\),
with \(\phi_v\) the \(N(0, v)\) density, is the Gaussian smoothing of
\(\Lam\) at time \(v\). One further Gaussian integration in \(\eta\) therefore
advances the time variable, and variances add along paths. Induction extends \eqref{eq:collapse} to \(k\)
sequential mediators at one certified evaluation in total. Adding a
mediator-mediator interaction \(\theta_{23} m_1 m_2\) makes the inner variance
\((\theta_3 + \theta_{23} m_1)^2 \sigma_2^2\) depend on \(m_1\). The time
parameter is then no longer constant across the convolution, so the collapse
degrades to \((k-1)\)-dimensional quadrature over the certified
one-dimensional core. For \(k\) parallel mediators with joint normal law
\(M = (M_1, \ldots, M_k)^{\top} \sim N(\mu, \Sigma)\) entering the outcome
model through \(a + b^{\top} M\), normality of linear combinations gives
\(\eta = a + b^{\top} \mu\) and \(v = b^{\top} \Sigma\, b\), so the certified
evaluation cost does not depend on \(k\). At every layer the uniform bound
survives composition, since the expectation of an error bounded by
\(\Delta_K^\star\) is itself bounded by \(\Delta_K^\star\), and every identity in
this subsection is an identity of the fixed closed form
\eqref{eq:int-from-mixture}.

\section{Conclusions}
\label{sec:conc}

We proved the Monahan--Stefanski equioscillation conjecture at every interior
minimiser of the normal-scale-mixture approximation of the logistic sigmoid,
through the Haar property of the mixture tangent space. An elementary
zero-counting difference argument then promoted every certified equioscillating
minimum to the unique global minimax. Interior attainment (\cref{prop:interior}),
and with it the conjecture, is established unconditionally for \(K \le 11\). We
certified the minimax constants for \(K = 9, 10, 11\), reaching a worst-case error
of \(1.26 \times 10^{-11}\). Beyond \(K = 11\) we certified mixtures with supremum error
strictly below \(\Delta_{11}^\star\) at \(K = 14, 15\), proved that the minimax
error vanishes in the limit (\cref{prop:limit}), and certified a strict signed
\(30\)-point equioscillation at \(K = 22\). For the
logistic-normal integral this yields a closed-form evaluation with a
parameter-free error bound, which in causal mediation needs no rare-outcome
assumption and supplants stochastic or grid-based integration at a cost of
milliseconds. The evaluation composes, so \(k\) sequential or parallel normal
mediators collapse to one certified evaluation of the same closed form
(\cref{sec:multi}). A conjecture that stood for thirty-four years is now a theorem
over the certified range, and the integral it concerns carries a certified closed
form within a factor of about thirty-three of double-precision unit roundoff. Four
questions remain open. The first is a proof or refutation of the
\(K = 12, 13\) plateau (\cref{conj:plateau}). The second is a strict
certificate at \(K > 11\), which appears to require global or topological
existence machinery rather than a local solver. The third is a rigorous
characterisation of the relaxed Bernstein--Akhiezer certificate frontier, a
question for classical approximation theory beyond the present scope. The fourth
is the general-\(K\) tangent-Wronskian identity (\cref{conj:wronskian}).

\section*{Declarations}

\paragraph*{Funding}
This work was supported by the National Medical Research Council, Singapore,
under grant MOH-SEEDFD24jul-0001.

\paragraph*{Competing interests}
The author has no competing interests to declare that are relevant to the
content of this article.

\paragraph*{Data availability}
No new data were collected. The public data set of \cref{sec:rare-outcome} is
distributed in the R package \texttt{faraway}, and the remaining computations
use synthetic data as described in \cref{sec:mediation}.

\paragraph*{Code availability}
The certified constants and the reproduction scripts are available at
\url{https://github.com/li-ruijie/LogisNormIntegral}, with a permanent DOI to be
minted before final acceptance.

\paragraph*{Author contributions}
R.L. is the sole author of this article and performed the theoretical analysis,
algorithm design, implementation, numerical experiments, and manuscript
preparation.

\begin{appendices}

\numberwithin{equation}{section}
\numberwithin{table}{section}
\numberwithin{figure}{section}

\crefalias{section}{appendix}
\crefalias{subsection}{appendix}
\crefname{appendix}{Appendix}{Appendices}
\Crefname{appendix}{Appendix}{Appendices}

\section{Setup: The family, the tangent space, and boundary values}
\label{sec:S-setup}

This appendix fixes the notation, the tangent space, and the boundary behaviour
of the error that every later proof uses. Fix \(K \ge 1\). The \(K\)-term
normal-scale-mixture family is
\begin{equation} \label{eq:S-family}
    \begin{split}
        \F_K = \Bigl\{\, f(\,\cdot\,; w, s) = \sum_{i=1}^K w_i\,\Phi(s_i\,\cdot\,)
        \;:\; {}& w_i > 0,\ \textstyle\sum_{i=1}^K w_i = 1,\\
        & 0 < s_1 < s_2 < \dots < s_K \,\Bigr\}.
    \end{split}
\end{equation}
Its closure \(\overline{\F_K}\) allows \(w_i \ge 0\) and coincident or boundary
scales \(s_i \in [0, \infty]\). An interior point of \(\F_K\) is one with all
\(w_i > 0\) and distinct positive finite scales. The constraint
\(\sum_i w_i = 1\) removes one weight degree of freedom, so \(\F_K\) is a smooth
manifold of dimension \(2K - 1\).

At an interior point \(f \in \F_K\) the tangent space is spanned by the partial
derivatives of \(f\) with respect to the free parameters,
\begin{equation} \label{eq:S-tangent}
    V = T_f\F_K = \operatorname{span}
    \bigl\{\, \Phi(s_1 z), \dots, \Phi(s_K z),\
    z\,\phi(s_1 z), \dots, z\,\phi(s_K z) \,\bigr\},
\end{equation}
where \(\partial f / \partial w_i = \Phi(s_i z)\) and \(\partial f / \partial s_i =
w_i\,z\,\phi(s_i z)\). The simplex constraint \(\sum_i w_i = 1\) imposes the single
linear relation \(\sum_i a_i = 0\) on the \(\Phi\)-coefficients, so
\(\dim V = 2K - 1\).

\begin{lemma}[Boundary values of tangent vectors] \label{lem:S-boundary}
    Every \(h \in V\), written \(h(z) = \sum_{i=1}^K a_i\,\Phi(s_i z) +
    \sum_{i=1}^K b_i\,z\,\phi(s_i z)\) with \(\sum_i a_i = 0\), satisfies
    \(h(0) = 0\) and \(\lim_{z \to \infty} h(z) = 0\).
\end{lemma}

\begin{proof}
    Since \(\Phi(0) = \tfrac12\) and \(0 \cdot \phi(0) = 0\), \(h(0) = \tfrac12
    \sum_i a_i = 0\). As \(z \to \infty\), \(\Phi(s_i z) \to 1\) and
    \(z\,\phi(s_i z) \to 0\), so \(\lim_{z \to \infty} h(z) = \sum_i a_i = 0\).
\end{proof}

These two boundary zeros are used in \cref{sec:S-haar,sec:S-uniqueness}, where they
convert an interior zero count of \(h\) into a zero count of \(h'\) through Rolle's
theorem.

\section{Existence and interior attainment of the minimax}
\label{sec:S-existence}

A point of the compactified parameter set with error at most \(\Delta_K^\star\) exists by
parameter compactness alone. The certified strict decrease then excludes the count-reducing limits, the
certified ratio bound excludes the limiting faces through \cref{lem:S-face}, and an induction
from the \(K = 1\) base places the minimiser in the open family for \(K \le 11\).

\begin{proposition}[Existence] \label{prop:S-existence}
    For every \(K \ge 1\) some \(f^\star \in \overline{\F_K}\) satisfies
    \(\sup_z|\Lam - f^\star| \le \Delta_K^\star\). For \(K = 1, \dots, 11\) the
    infimum is moreover attained in \(\F_K\) itself. Some interior \(f^\star\) has
    \(\sup_z|\Lam - f^\star| = \Delta_K^\star\), every minimiser lies in \(\F_K\), and every
    point of \(\overline{\F_K}\) has uniform error at least \(\Delta_K^\star\).
\end{proposition}

\begin{proof}
    Let \((f_n) \subset \F_K\) be a minimising sequence, \(f_n = \sum_i w_i^{(n)}
    \Phi(s_i^{(n)} z)\), with \(\sup_z|\Lam - f_n| \to \Delta_K^\star\). The weights
    lie in the compact simplex, so after passing to a subsequence \(w_i^{(n)} \to
    w_i^\star \in [0, 1]\) with \(\sum_i w_i^\star = 1\), and each scale sequence
    lies in the compact \([0, \infty]\), so a further subsequence gives
    \(s_i^{(n)} \to s_i^\star \in [0, \infty]\). Parameter convergence forces
    \(f_n \to f^\star = \sum_i w_i^\star \Phi(s_i^\star z) \in \overline{\F_K}\)
    pointwise at every \(z \ne 0\), under the conventions \(\Phi(0 \cdot z) =
    \tfrac12\) and \(\Phi(\infty \cdot z) = \mathbf{1}\{z > 0\}\). Each summand
    is continuous in its parameters except at the single point \(z = 0\) of
    the step limit. Hence \(|\Lam(z) - f^\star(z)| = \lim_n |\Lam(z) - f_n(z)| \le
    \Delta_K^\star\) for every \(z \ne 0\), and at the origin \(|\Lam(0) -
    f^\star(0)| = m_\infty/2 = \lim_{z \downarrow 0}|\Lam(z) - f^\star(z)| \le
    \Delta_K^\star\), where \(m_\infty\) is the limit mass at \(s = \infty\).
    Therefore \(\sup_z|\Lam - f^\star| \le \Delta_K^\star\).

    Whether this supremum
    \emph{equals} \(\Delta_K^\star\) is a separate question. A priori
    \(\overline{\F_K}\) could contain points of error strictly below the infimum
    over the open family, since uniform approximation of a step by continuous
    mixtures fails near the origin, and only the certified constants close the
    gap. For \(K \le 11\) they do, by induction on \(K\). Assume every point of
    \(\overline{\F_{K-1}}\) has error at least \(\Delta_{K-1}^\star\) (for
    \(K = 1\) this is vacuous). A count-reducing limit, a vanishing weight or
    coincident scales, lies in \(\overline{\F_{K-1}}\), so for \(K \ge 2\) its error
    is at least \(\Delta_{K-1}^\star > \Delta_K^\star\) by the certified strict
    decrease (\cref{tab:cert}), while at \(K = 1\) no such limit exists, since
    \(\sum_i w_i = 1\) forces the single weight to one.

    Boundary mass at \(s \in \{0, \infty\}\) is excluded for
    \(K = 1\) directly, since the only such points are the constant \(\tfrac12\) and
    the unit step, each with error \(\tfrac12\), and for \(K \ge 2\) by
    \cref{prop:S-interior,lem:S-face}, whose closure hypothesis is the same
    inductive bound at level \(K - 1\). Hence every point of
    \(\overline{\F_K} \setminus \F_K\) has error strictly above \(\Delta_K^\star\),
    the limit \(f^\star\) lies in \(\F_K\) at the value \(\Delta_K^\star\), and
    every point of \(\overline{\F_K}\) has error at least \(\Delta_K^\star\), which
    advances the induction.
\end{proof}

The equioscillation theorem is stated at an \emph{interior} minimiser. Interior
attainment is governed by strict monotonicity of the minimax sequence.

\begin{proposition}[Interior attainment under strict monotonicity]
\label{prop:S-interior}
    Fix \(K \ge 2\) with \(\Delta_K^\star < \Delta_{K-1}^\star\), and suppose every point of
    \(\overline{\F_{K-1}}\) has uniform error at least \(\Delta_{K-1}^\star\). Then no minimiser
    \(f^\star\) has a vanishing weight or two coincident scales. If in addition the
    minimax error over each limiting face (a constant, a unit step, or both, atop
    at most \(K - 1\) genuine Gaussians) strictly exceeds \(\Delta_K^\star\), then
    every minimiser is interior. All three hypotheses hold for \(K = 2, \dots, 11\): the
    strict inequality from the certified constants (\cref{tab:cert}), the closure bound from
    the induction of \cref{prop:S-existence}, and the face condition from \cref{lem:S-face}.
\end{proposition}

\begin{proof}
    A count-reducing degeneracy, a vanishing weight or two coincident finite
    positive scales, leaves a mixture of at most \(K - 1\) distinct Gaussians, so
    \(f^\star \in \overline{\F_{K-1}}\) and its error is at least
    \(\Delta_{K-1}^\star > \Delta_K^\star\), by the closure hypothesis and the strict
    decrease. Such a point does not attain the minimum, which proves the first claim.

    The limiting degeneracies place mass on the scales \(s = 0\), where \(\Phi(0
    \cdot z) \equiv \tfrac12\) contributes a constant, and \(s = \infty\), where
    \(\Phi(\infty \cdot z) = \mathbf 1\{z > 0\}\) contributes a unit step at the
    origin. Write \(m_0\) and \(m_\infty\) for the mass on these two scales. As
    \(z \to +\infty\) the finite atoms saturate to their weights, the \(s = 0\)
    atoms sit at half of theirs, and the \(s = \infty\) atoms reach theirs, giving
    \(f^\star(+\infty) = 1 - m_0/2\) against \(\Lam(+\infty) = 1\), a tail error of
    \(m_0/2\). At \(z = 0\) the step atoms jump by \(m_\infty\) across the
    continuous value \(\Lam(0) = \tfrac12\), forcing an error of at least
    \(m_\infty/2\). Hence \(\sup_z|\Lam - f^\star| \ge \tfrac12\max(m_0,
    m_\infty)\), so any limiting-face minimiser carries vanishing mass
    \(\max(m_0, m_\infty) \le 2\Delta_K^\star\) atop at most \(K - 1\) genuine
    Gaussians. The hypothesis that each limiting face exceeds \(\Delta_K^\star\)
    excludes this residual, and every minimiser is then interior.
\end{proof}

\begin{lemma}[Limiting faces exceed the interior value] \label{lem:S-face}
    Fix \(K \ge 2\) with \(\Delta_K^\star / \Delta_{K-1}^\star < (1 -
    4\Delta_K^\star)/2\), and suppose every point of \(\overline{\F_{K-1}}\) has
    uniform error at least \(\Delta_{K-1}^\star\). Then the minimax error over
    every limiting face of \(\overline{\F_K}\) strictly exceeds
    \(\Delta_K^\star\).
\end{lemma}

\begin{proof}
    Let \(f\) be any point of a limiting face with \(\|\Lam - f\| \le
    \Delta_K^\star\), and write \(m_0\), \(m_\infty\) for its mass at \(s = 0\) and
    \(s = \infty\). The tail bound \(\sup_z|\Lam - f| \ge \tfrac12\max(m_0,
    m_\infty)\) gives \(\max(m_0, m_\infty) \le 2\Delta_K^\star\), so the limiting
    mass \(\mu = m_0 + m_\infty \le 4\Delta_K^\star\). Renormalising the at most
    \(K - 1\) genuine Gaussians yields \(g = (f - \tfrac12 m_0 - m_\infty
    \mathbf 1\{z > 0\})/(1 - \mu) \in \overline{\F_{K-1}}\), so \(\|\Lam - g\| \ge
    \Delta_{K-1}^\star\) by the closure hypothesis. With
    \(h = \tfrac12 m_0 + m_\infty \mathbf 1\{z > 0\} -
    \mu\Lam\), a term-by-term check on \(z \le 0\) and \(z > 0\) gives \(\|h\| \le
    \tfrac12\max(m_0, m_\infty) \le \Delta_K^\star\), and the identity
    \((1 - \mu)(\Lam - g) = (\Lam - f) + h\) yields
    \[
        \|\Lam - f\| \ge (1 - \mu)\,\Delta_{K-1}^\star - \Delta_K^\star
        \ge (1 - 4\Delta_K^\star)\,\Delta_{K-1}^\star - \Delta_K^\star .
    \]
    The hypothesis makes the right-hand side strictly exceed \(\Delta_K^\star\),
    contradicting \(\|\Lam - f\| \le \Delta_K^\star\).
\end{proof}

For \(K = 2, \dots, 11\) the certified solutions satisfy \(\Delta_K^\star <
\Delta_{K-1}^\star\) strictly, giving the strict-decrease hypothesis of
\cref{prop:S-interior}. The certified ratios
\(\Delta_K^\star/\Delta_{K-1}^\star\) are at most \(0.19\) throughout \(2 \le K
\le 11\), whereas the threshold \((1 - 4\Delta_K^\star)/2\) exceeds \(0.498\) (since
\(\Delta_K^\star \le \Delta_2^\star < 6 \times 10^{-4}\)), so the ratio hypothesis
of \cref{lem:S-face} holds at every level.

Its closure hypothesis, shared with \cref{prop:S-interior}, is discharged by
the induction of \cref{prop:S-existence}, climbing from the \(K = 1\) base (the
constant \(\tfrac12\) and the unit step, each with error \(\tfrac12 >
\Delta_1^\star\)) one level at a time. At each level the lemma excludes the faces,
the strict decrease excludes the count-reducing limits, and the closure bound
emerges for the next level's use. All three hypotheses of \cref{prop:S-interior}
therefore hold for \(K = 2, \dots, 11\), and with the \(K = 1\) case settled
directly in the proof of \cref{prop:S-existence}, the minimiser is interior for
\(K = 1, \dots, 11\). Strictness
and the ratio condition beyond \(K = 11\) are open, and the equioscillation
theorem below is accordingly conditioned on interior attainment.

\section{The Haar property of the tangent space}
\label{sec:S-haar}

The zero-counting engine is the following classical bound on the real zeros of a
generalised exponential sum.

\begin{lemma}[Zeros of generalised exponential sums] \label{lem:S-exp-sum-zeros}
    Let \(g\) be the generalised exponential sum
    \(g(x) = \sum_{i=1}^m P_i(x)\,e^{\lambda_i x}\), with distinct real
    \(\lambda_i\) and polynomials \(P_i \not\equiv 0\) of degree \(r_i\). Then \(g\)
    has at most \(\bigl(\sum_{i=1}^m (r_i + 1)\bigr) - 1\) real zeros, counted with
    multiplicity. In particular a plain exponential sum of \(N\) terms has at most
    \(N - 1\) zeros, and a multiplicity-two sum \(\sum_{i=1}^K (c_i + d_i x)\,e^{\lambda_i x}\)
    at most \(2K - 1\).
\end{lemma}

\begin{proof}
    These bounds are classical~\cite{polya-ueberan-1913, schoenbe-polyafr-1951,
    karlin-tchebych-1966}. Induct on the total dimension \(M = \sum_i (r_i + 1)\).
    If \(M = 1\) then \(g = c\,e^{\lambda x}\) is nowhere zero. For \(M > 1\),
    multiplying by \(e^{-\lambda_1 x}\) leaves the zeros unchanged, so take
    \(\lambda_1 = 0\) and write \(g = P_1 + \sum_{i \ge 2} P_i\,e^{\lambda_i x}\).
    If \(m = 1\) then \(g = P_1\) is a polynomial of degree \(r_1 = M - 1\) and the
    bound is immediate, so assume \(m \ge 2\).
    Differentiating \(r_1 + 1\) times annihilates \(P_1\), while each remaining term
    maps to \(\tilde P_i\,e^{\lambda_i x}\) with \(\deg \tilde P_i = r_i\). Thus
    \(g^{(r_1 + 1)}\) is a sum of the same form with total dimension
    \(M - (r_1 + 1)\), hence has at most \(M - (r_1 + 1) - 1\) zeros by induction.
    Applying Rolle's theorem \(r_1 + 1\) times gives at most \(M - 1\) zeros for
    \(g\).
\end{proof}

\begin{theorem}[Haar property of the tangent space] \label{thm:S-haar}
    The tangent space \(V\) of~\eqref{eq:S-tangent} is a Chebyshev (Haar) system of
    dimension \(2K - 1\) on \((0, \infty)\), that is, every nonzero \(h \in V\) has at most
    \(2K - 2\) distinct zeros in \((0, \infty)\).
\end{theorem}

\begin{proof}
    Let \(h \in V\) be nonzero, written as in \cref{lem:S-boundary}, so \(h(0) =
    0\) and \(h(z) \to 0\) as \(z \to \infty\). Suppose, for contradiction, that
    \(h\) has \(N \ge 2K - 1\) distinct zeros in \((0, \infty)\). Adjoining the
    boundary zeros at \(0\) and \(\infty\), Rolle's theorem gives
    \(\#\{z \in (0, \infty) : h'(z) = 0\} \ge N + 1\). Differentiating and using
    \(\phi'(t) = -t\,\phi(t)\),
    \[
        h'(z) = \sum_{i=1}^K (a_i s_i + b_i)\,\phi(s_i z)
              - \sum_{i=1}^K b_i s_i^2\,z^2\,\phi(s_i z).
    \]
    Substituting \(\phi(s_i z) = (2\pi)^{-1/2} e^{-s_i^2 z^2 / 2}\), setting
    \(x = z^2\) and \(\lambda_i = s_i^2 / 2 > 0\) (all distinct), zeros of \(h'\) on
    \((0, \infty)\) biject with zeros on \((0, \infty)\) of
    \[
        g(x) = \sum_{i=1}^K c_i\,e^{-\lambda_i x}
             + \sum_{i=1}^K d_i\,x\,e^{-\lambda_i x}
    \]
    for scalars \(c_i, d_i\) depending linearly on \((a_i, b_i)\). The scalars
    are not all zero, since the map \((a, b) \mapsto (c, d)\) is invertible for
    positive scales, so \(h \ne 0\) forces \(g \ne 0\). By the multiplicity-two case of
    \cref{lem:S-exp-sum-zeros}, every nonzero \(g\) of this form has at most
    \(2K - 1\) zeros. Equivalently, the family
    \(\{e^{-\lambda_i x},\, x\,e^{-\lambda_i x}\}_{i=1}^K\) is a \(2K\)-dimensional
    Extended Complete Chebyshev system on \([0, \infty)\)~\cite{karlin-totalpo-1968,
    braess-nonlinea-1986}. Therefore \(N + 1 \le 2K - 1\), that is \(N \le 2K - 2\),
    contradicting \(N \ge 2K - 1\).
\end{proof}

\section{Lower bound: At least \texorpdfstring{\(2K\)}{2K} alternations}
\label{sec:S-lower}

The Kolmogorov variational criterion, combined with the Haar property proved in
\cref{sec:S-haar}, gives a lower bound on the number of alternations that an interior
minimiser must carry.

\begin{corollary}[Equioscillation, lower bound] \label{cor:S-equiosc}
    Let \(f^\star \in \F_K\) attain the minimax \(\Delta_K^\star\) at an interior
    point (\cref{prop:S-existence,prop:S-interior}). Then the error
    \(e^\star = \Lam - f^\star\) has
    at least \(2K\) alternating extrema in \((0, \infty)\), that is, there exist
    \(0 < z_1 < \dots < z_{2K}\) with \(|e^\star(z_j)| = \Delta_K^\star\) and
    \(\sgn\,e^\star(z_j) = (-1)^{j+1}\sgn\,e^\star(z_1)\).
\end{corollary}

\begin{proof}
    Note first that \(\Delta_K^\star > 0\). A zero error would make \(\Lam\) a
    finite normal-scale mixture, impossible since \(1 - \Lam(z)\) decays like
    \(e^{-z}\) while every mixture tail decays at the Gaussian rate of its
    smallest scale. Since \(\Lam\) and the Gaussians are real-analytic and
    \(e^\star \to 0\) at \(+\infty\), the extremal set \(E = \{z \in (0, \infty) :
    |e^\star(z)| = \Delta_K^\star\}\) is finite. Order its points and group them
    into maximal runs of constant \(\sgn\,e^\star\), and let \(m\) be the number
    of runs. Suppose \(m \le 2K - 1\), and pick a separating knot between each
    adjacent pair of oppositely signed runs, giving \(t_1 < \dots < t_{m-1}\) in
    gaps where \(|e^\star| < \Delta_K^\star\). By \cref{thm:S-haar}, \(V\) is a
    Haar space of dimension \(2K - 1 \ge m\), so by the standard Chebyshev
    construction it contains an \(h\) whose only sign changes are at the
    \(m - 1 \le 2K - 2\) knots, with \(\sgn\,h \equiv \sgn\,e^\star\) on each
    run~\cite{rice-approxim-1969, cheney-introduc-1982}. Then \(e^\star\,h > 0\)
    on all of \(E\).

    The tangent direction is realised by a parameter move, which is where
    interiority enters. Writing \(h(z) = \sum_i a_i \Phi(s_i^\star z) + \sum_i
    b_i\, z\,\phi(s_i^\star z)\) with \(\sum_i a_i = 0\), the curve
    \(f_\varepsilon = f(\,\cdot\,;\, w_i^\star + \varepsilon a_i,\ s_i^\star +
    \varepsilon b_i / w_i^\star)\) stays in \(\F_K\) for all small \(\varepsilon
    \ge 0\), since the weights stay positive with unit sum and the scales stay
    ordered, as the family equation~\eqref{eq:S-family} requires. The curve
    satisfies \(f_\varepsilon = f^\star + \varepsilon h + R_\varepsilon\) with
    \(\sup_{z \ge 0} |R_\varepsilon(z)| = O(\varepsilon^2)\).
    The remainder is uniform in \(z\) because every second parameter
    derivative of the mixture is bounded on \([0, \infty)\).

    Fix a closed
    neighbourhood \(U\) of \(E\) on which \(e^\star h \ge c > 0\). On \(U\),
    \(|e^\star - \varepsilon h|^2 \le (\Delta_K^\star)^2 - 2 \varepsilon c +
    O(\varepsilon^2)\), so \(|\Lam - f_\varepsilon| \le \Delta_K^\star -
    \varepsilon c / \Delta_K^\star + O(\varepsilon^2) < \Delta_K^\star\) for small
    \(\varepsilon > 0\). Off \(U\), the continuous \(|e^\star|\) stays below
    \(\Delta_K^\star - \delta_0\) for some \(\delta_0 > 0\), since its extremal set
    lies inside \(U\) and \(e^\star\) vanishes at infinity, and the
    \(O(\varepsilon)\) move preserves the gap. Hence \(\sup_z|\Lam -
    f_\varepsilon| < \Delta_K^\star\) for all sufficiently small \(\varepsilon >
    0\), contradicting minimaxity. This parameter-space perturbation is the
    Kolmogorov criterion for a tangentially Haar
    manifold~\cite{rice-approxim-1969, cheney-introduc-1982}, carried out
    explicitly. Hence \(m \ge 2K\), and one extremum from each run gives
    \(0 < z_1 < \dots < z_{2K}\) of magnitude \(\Delta_K^\star\) with strictly
    alternating signs.
\end{proof}

\section{Upper bound: At most \texorpdfstring{\(2K\)}{2K} extrema, via total
    positivity}
\label{sec:S-upper}

This section establishes the matching upper bound of \(2K\) extrema directly, by
counting the zeros of the differentiated error as a mixed Laplace transform.
Neither positivity nor optimality of the weights is used, so the bound holds for
arbitrary signed mixtures.

\begin{theorem}[Equioscillation, upper bound and exact count] \label{thm:S-upper}
    Let \(f(z) = \sum_{i=1}^K w_i \Phi(s_i z)\) with real weights \(w_i\) and
    distinct positive scales \(s_i\). Then \(e = \Lam - f\) has at most \(2K\)
    local extrema in \((0, \infty)\). In particular, under the hypotheses of
    \cref{cor:S-equiosc} the minimax error \(e^\star\) has at most \(2K\) interior
    extrema, so combined with \cref{cor:S-equiosc} it equioscillates at exactly
    \(2K\) points of \((0, \infty)\), and (by \cref{sec:S-symmetry}) at \(4K\)
    points of \(\R\). This count resolves the Monahan--Stefanski conjecture for the
    normal-scale-mixture family at every interior minimiser, hence unconditionally
    for \(K \le 11\) (\cref{prop:S-existence}).
\end{theorem}

\begin{proof}
    By the exact normal-scale-mixture representation~\cite{andrews-scalemi-1974, stefansk-normals-1991} of
    the logistic sigmoid~\cite[eq.~18.4.6, written there in the scale
    \(\sigma = 1/s\)]{monahan-normals-1992},
    \begin{equation} \label{eq:S-stefanski}
        \Lam(z) = \int_0^\infty \Phi(s z)\,q(s)\,ds,
        \qquad q > 0,\quad \int_0^\infty q(s)\,ds = 1,
    \end{equation}
    with \(q\) the Kolmogorov--Smirnov mixing density (the Stefanski mixing density
    of the literature). The bound
    \(\int_0^\infty s\,q(s)\,ds < \infty\) justifies differentiation under the
    integral sign and, for \(z > 0\),
    \begin{equation} \label{eq:S-deriv}
        e'(z) = \int_0^\infty s\,\phi(s z)\,q(s)\,ds
        - \sum_{i=1}^K w_i s_i\,\phi(s_i z).
    \end{equation}
    Under \(u = z^2/2\) and \(\lambda = s^2\) the Gaussian density becomes the
    Laplace kernel, \(\phi(s z) = (2\pi)^{-1/2} e^{-\lambda u}\), and
    \(\sqrt{2\pi}\, e'\) becomes the mixed Laplace transform
    \begin{equation} \label{eq:S-transform}
        F(u) = \int_0^\infty \rho(\lambda)\, e^{-\lambda u}\,d\lambda
        - \sum_{i=1}^K c_i\, e^{-\lambda_i u},
        \qquad \rho(\lambda) = \tfrac12\, q(\sqrt{\lambda}\,),
        \quad c_i = w_i s_i,
    \end{equation}
    with \(\lambda_i = s_i^2\) and \(\rho > 0\). The polynomial moments of
    \(\rho\) against
    every \(e^{-\lambda u}\), \(u > 0\), are finite, since the exponential factor
    absorbs them. By \cref{rem:S-peel} below,
    \(F\) has at most \(2K\) zeros in \((0, \infty)\), whatever the signs of the
    \(c_i\). This bound is the zero-counting form of the variation-diminishing property
    of the sign-regular Laplace kernel~\cite[Chapters~I
    and~V]{karlin-totalpo-1968}. Since \(u = z^2/2\) is a bijection of
    \((0, \infty)\), \(e'\) has at most \(2K\) zeros there, and every interior
    local extremum of \(e\) is a zero of \(e'\), so \(e\) has at most \(2K\) local
    extrema in \((0, \infty)\).

    At an interior minimiser, every point of the extremal set \(\{z > 0 :
    |e^\star(z)| = \Delta_K^\star\}\) is a global, hence local, extremum of
    \(e^\star\), so that set has at most \(2K\) points. With the lower bound of
    \cref{cor:S-equiosc}, the count equals \(2K\).
\end{proof}

\begin{remark}[The double-peel zero count] \label{rem:S-peel}
    The bound on the mixed Laplace transform~\eqref{eq:S-transform} is the Rolle
    engine of \cref{lem:S-exp-sum-zeros} lifted from finite sums to
    density-plus-atoms transforms. Let \(F(u) = \int_0^\infty
    \rho(\lambda)\, e^{-\lambda u}\,d\lambda + \sum_{i=1}^K c_i\,
    e^{-\lambda_i u}\) with an almost-everywhere positive density and \(K\) atoms
    of arbitrary sign. Multiplying by \(e^{\lambda_1 u}\) leaves the zeros
    unchanged, and differentiating costs at most one zero by Rolle's theorem
    while replacing \((\rho(\lambda), c_i)\) by \((\rho(\lambda)(\lambda_1 - \lambda),\,
    c_i(\lambda_1 - \lambda_i))\). The first atom is thus annihilated, at the
    price that the density factor now changes sign at \(\lambda_1\). A second identical
    peel multiplies the density by \((\lambda_1 - \lambda)\) again, restoring the
    almost-everywhere positive density \(\rho(\lambda)(\lambda_1 - \lambda)^2\)
    and keeping the first atom dead. Two Rolle steps thus remove one atom, so the
    \(K\) atoms cost \(2K\) zeros, and the atom-free base case is a strictly
    positive integral with no zeros. Hence \(F\) has at most \(2K\) zeros in
    \((0, \infty)\).
\end{remark}

\section{Uniqueness and the tangent-Wronskian identity}
\label{sec:S-uniqueness}

Uniqueness is elementary and self-contained, resting only on the zero-counting
engine of \cref{sec:S-haar} applied to the difference of two competing mixtures.
In the classical vocabulary the same conclusion follows from varisolvency, since
a full-dimension Haar tangent at every interior point makes the family
varisolvent of degree \(2K - 1\). For a varisolvent family the Kolmogorov
criterion is necessary and sufficient and best approximations are
unique~\cite[Theorem~7-4]{rice-approxim-1969},
\cite[III.4.6]{braess-nonlinea-1986}. The direct argument below keeps these
appendices free of that machinery.

\begin{lemma}[Difference zero bound] \label{lem:S-diff}
    Let \(f, g \in \F_K\) be distinct. Then \(d = f - g\) is not identically zero
    and has at most \(2K - 2\) zeros in \((0, \infty)\), counted with
    multiplicity.
\end{lemma}

\begin{proof}
    Merging coincident scales across the two parameter lists, \(d(z) = \sum_j
    c_j\,\Phi(\sigma_j z)\) is a combination of \(m \le 2K\) distinct
    \(\Phi\)-atoms with \(\sum_j c_j = 1 - 1 = 0\). Not every \(c_j\) vanishes,
    since all coefficients vanishing would force the two mixtures to share
    scales and weights. Hence \(d \not\equiv 0\), for \(d \equiv 0\) would force
    \(d' \equiv 0\), while \(d'(z) = \sum_j c_j \sigma_j\,\phi(\sigma_j z)\) maps
    under \(x = z^2\) to a nontrivial plain exponential sum, which by
    \cref{lem:S-exp-sum-zeros} has finitely many zeros. As in \cref{lem:S-boundary}, \(d(0) = 0\) and
    \(d(z) \to 0\) as \(z \to \infty\). If \(d\) had \(2K - 1\) zeros in
    \((0, \infty)\) with multiplicity, Rolle's theorem across these and the two
    boundary zeros would give \(d'\) at least \(2K\) zeros in \((0, \infty)\) with
    multiplicity. The exponential-sum image of \(d'\) has \(m \le 2K\)
    terms and so, by \cref{lem:S-exp-sum-zeros}, at most \(m - 1 \le 2K - 1\)
    zeros with multiplicity.
\end{proof}

\begin{proposition}[de la Vall\'ee Poussin bound and its equality case]
\label{prop:S-dlvp}
    Let \(f \in \F_K\) have error \(e_f = \Lam - f\) alternating at \(2K\) points
    \(0 < z_1 < \dots < z_{2K}\), in the sense that for some \(\delta > 0\) and
    \(\epsilon = \pm 1\), \((-1)^{j} \epsilon\, e_f(z_j) \ge \delta\) for every
    \(j\). Then \(\Delta_K^\star \ge \delta\), and every \(g \in \F_K\) with
    \(\sup_{z}|\Lam - g| \le \delta\) equals \(f\).
\end{proposition}

\begin{proof}
    Let \(g \in \F_K\) satisfy \(\sup_z|\Lam - g| \le \delta\) and suppose
    \(g \ne f\). At each alternation point, \(d(z_j) = (f - g)(z_j) = e_g(z_j) -
    e_f(z_j)\) with \(|e_g(z_j)| \le \delta \le (-1)^{j} \epsilon\, e_f(z_j)\),
    so \((-1)^{j} \epsilon\, d(z_j) \le 0\), and the signs of \(d\) weakly
    alternate against the \(2K\) points. A function that is not identically zero,
    is analytic, and whose signs weakly alternate against \(2K\) ordered points
    has at least \(2K - 1\) zeros in \([z_1, z_{2K}]\) counted with multiplicity.
    Each consecutive pair contributes a sign-change zero or an endpoint zero. An
    endpoint zero claimed by both neighbouring pairs is either a sign crossing,
    in which case one of the two adjacent open intervals holds a further zero, or
    a zero without sign change, which has even multiplicity and counts twice.
    This count contradicts the bound \(2K - 2\) of \cref{lem:S-diff}, so \(g = f\).
    Consequently \(\Delta_K^\star \ge \delta\), since a mixture \(g\) with
    \(\sup_z|\Lam - g| < \delta\) would equal \(f\) by the above while the error
    of \(f\) reaches \(\delta\) in modulus.
\end{proof}

\begin{corollary}[Global uniqueness] \label{cor:S-global}
    At every \(K\) for which the minimax is attained at an interior point (in
    particular \(K = 1, \dots, 11\), by \cref{prop:S-existence}), the minimiser of
    \(\Delta_K^\star\) is unique in \(\F_K\), and any interior mixture whose
    error alternates at \(2K\) points at level \(\Delta_K^\star\) is that
    minimiser.
\end{corollary}

\begin{proof}
    An interior minimiser \(f^\star\) equioscillates at \(2K\) points at level
    \(\Delta_K^\star\) (\cref{cor:S-equiosc}). Any other minimiser \(g\) has
    \(\sup_z|\Lam - g| = \Delta_K^\star\), so \(g = f^\star\) by the equality
    case of \cref{prop:S-dlvp} with \(\delta = \Delta_K^\star\).
\end{proof}

\begin{corollary}[Strict local uniqueness] \label{cor:S-unique}
    Under the hypotheses of \cref{cor:S-global} the minimiser \((w^\star,
    s^\star)\) is strict, in the sense that every other mixture in \(\F_K\) has
    uniform error strictly above \(\Delta_K^\star\). In particular it is a strict, locally
    isolated minimum of the minimax objective, the conclusion the classical
    Remez-Jacobian route reaches through the implicit function
    theorem~\cite{rice-approxim-1969, cheney-introduc-1982}.
\end{corollary}

\begin{proof}
    By \cref{cor:S-global} the minimum is attained only at \(f^\star\), so every
    \(g \ne f^\star\) has \(\sup_z|\Lam - g| > \Delta_K^\star\).
\end{proof}

A separate explicit certificate for the strict local uniqueness, independent of
the equioscillation argument above and available at finite \(K\), is a structural
identity for the confluent tangent Wronskian. Let \(Q(z; s)\) denote
the reduced \((2K - 1) \times (2K - 1)\) tangent-Wronskian cofactor of the basis
\(\{\Phi(s_i z),\, z\,\phi(s_i z)\}_{i=1}^K\) after deleting the order-\(0\)
derivative row and the largest-index column \(z\,\phi(s_K z)\). Rows are ordered
by ascending derivative order \(1, \dots, 2K - 1\), and the surviving columns are
ordered as the \(\Phi\) block \(\Phi(s_1 z), \dots, \Phi(s_K z)\) followed by the
\(z\,\phi\) block \(z\,\phi(s_1 z), \dots, z\,\phi(s_{K-1} z)\). The sign of the
closed form below is tied to this ordering. Permuting rows or columns changes
\(Q\) only by a global sign, so the nonvanishing conclusion is
ordering-independent.

\begin{theorem}[Structural pure-monomial identity, \(K \le 5\)]
\label{thm:S-wronskian}
    For \(2 \le K \le 5\) and ordered positive distinct
    \(s_1 < s_2 < \dots < s_K\), the reduced tangent-Wronskian cofactor factors as
    \(Q(z; s) = c_K(s)\, z^{(2K-1)(K-1)}\), where \(c_K(s)\) is a signed product of a
    power of two, powers of the scales, and even powers of the pairwise differences
    of squared scales, nonzero on the ordered positive cone, so \(Q(z; s) \ne 0\) on
    \((0, \infty)\).
\end{theorem}

\begin{proof}
    For each \(K \in \{2, 3, 4, 5\}\) the reduced \((2K-1) \times (2K-1)\) tangent
    matrix has entries in \(\mathbb{Q}[z, s_1, \dots, s_K]\). A Bareiss fraction-free
    elimination~\cite{bareiss-sylveste-1968} (\(K \le 4\)) and a staged suffix-minor
    expansion in exact integer arithmetic with modular cross-checks at multiple
    points and primes (\(K = 5\)) return the determinant as an exact multivariate
    polynomial. The polynomial factors as
    \(c_K(s)\,z^{(2K-1)(K-1)}\) with
    \[
        c_K(s) = (-1)^{\lfloor (K-1)/2 \rfloor}\,2^{K-1}\,s_K
        \prod_{i=1}^{K-1} s_i^3
        \prod_{1 \le i < j \le K-1} (s_i^2 - s_j^2)^4
        \prod_{i=1}^{K-1} (s_i^2 - s_K^2)^2 .
    \]
    Each factor is nonzero for ordered positive distinct \(s\), so \(c_K(s) \ne 0\)
    and \(Q(z; s) \ne 0\) on \((0, \infty)\).
\end{proof}

The exponents in \(c_K(s)\) follow the confluent-Vandermonde count for the nodes
\(\lambda_i = s_i^2\). Each scale \(i < K\) contributes a block of size two (its
two tangent directions \(\Phi(s_i z)\) and \(z\,\phi(s_i z)\)), the
deleted-column scale \(s_K\) a block of size one. Each pair of nodes enters
with exponent equal to the product of its block sizes, \(2 \cdot 2 = 4\) within
the first \(K - 1\) scales and \(2 \cdot 1 = 2\) against \(s_K\). The conjectured
general-\(K\) form extrapolates this count. Its open content is the pure monomial
factorisation in \(z\).

\begin{conjecture}[Pure-monomial identity for general \(K\)]
\label{conj:S-wronskian}
    The factorisation of \cref{thm:S-wronskian} holds, with the same closed form
    for \(c_K(s)\), for every \(K \ge 2\).
\end{conjecture}

For \(6 \le K \le 12\) the formula is verified pointwise at the integer scales
\(s = (1, 2, \dots, K)\). The Bareiss determinant and the formula evaluate to the
same integer, which reaches \(140\) digits at \(K = 8\) and
\(412\) digits at \(K = 12\). A general-\(K\) symbolic proof is open but is not
required here. The uniqueness chain above rests only on the zero-counting engine
of \cref{sec:S-haar}, and at each certified \(s^\star\) the nonvanishing of
\(Q\) is in addition confirmed directly by a Sturm root count of \(Q\) over
\(\mathbb{Q}\).

\section{Symmetry: From \texorpdfstring{\((0, \infty)\)}{(0,infinity)} to
    \texorpdfstring{\(\R\)}{R}}
\label{sec:S-symmetry}

The counts established so far are on \((0, \infty)\). Oddness of the error transfers
them to the whole line, doubling every count.

\begin{lemma}[The minimax error is odd] \label{lem:S-odd}
    For any \(f \in \F_K\) the error \(e = \Lam - f\) satisfies \(e(-z) = -e(z)\)
    for all \(z \in \R\).
\end{lemma}

\begin{proof}
    The logistic sigmoid satisfies \(\Lam(-z) = 1 - \Lam(z)\), and \(\Phi(-u) = 1 -
    \Phi(u)\), so \(\Phi(-s_i z) = 1 - \Phi(s_i z)\). Using \(\sum_i w_i = 1\),
    \(f(-z) = \sum_i w_i(1 - \Phi(s_i z)) = 1 - f(z)\). Therefore \(e(-z) = (1 -
    \Lam(z)) - (1 - f(z)) = -\,e(z)\).
\end{proof}

The oddness of \(e^\star\) forces \(e^\star(0) = 0\), and the extrema occur in
mirror pairs. The \(2K\) equioscillation points of \cref{thm:S-upper} on
\((0, \infty)\) reflect to \(2K\) points on \((-\infty, 0)\), giving exactly
\(4K\) equioscillation points of the extended error on \(\R\), with no extremum at
the origin.

\section{An alternative (obstructed) total-positivity route}
\label{sec:S-alt-tp}

A classical alternative would read the equioscillation count off total-positivity
theory applied directly to the cumulative kernel \(\Phi(sz)\), bypassing the
differentiation of \cref{sec:S-upper}. This section records why that route is
obstructed. The cumulative kernel is log-concave but fails to be totally positive.

\begin{lemma}[Log-concavity of the normal distribution function] \label{lem:S-logconcave}
    The map \(z \mapsto \log\Phi(z)\) is strictly concave on \(\R\), and the
    reversed hazard rate \(z \mapsto \phi(z)/\Phi(z)\) is strictly decreasing.
\end{lemma}

\begin{proof}
    With \(\phi' = -z\,\phi\), \((\log\Phi)''(z) = -\phi(z)[z\,\Phi(z) +
    \phi(z)]/\Phi(z)^2\), so the sign equals that of \(-[z\,\Phi(z) + \phi(z)]\).
    Set \(g(z) = z\,\Phi(z) + \phi(z)\). Then \(g'(z) = \Phi(z) > 0\) and \(g(z) \to
    0\) as \(z \to -\infty\), so \(g > 0\), giving \((\log\Phi)'' < 0\). The same
    inequality makes \((\phi/\Phi)' = -g\,\phi/\Phi^2 < 0\).
\end{proof}

\begin{proposition}[The cumulative kernel is not totally positive]
\label{prop:S-not-tp}
    The kernel \(\Phi(sz)\) is not totally positive of order \(2\) on
    \((0, \infty) \times \R\). For \(0 < s_1 < s_2\) the minor \(D(z_1, z_2) =
    \Phi(s_1 z_1)\,\Phi(s_2 z_2) - \Phi(s_1 z_2)\,\Phi(s_2 z_1)\) is strictly
    negative at some ordered arguments and strictly positive at others.
\end{proposition}

\begin{proof}
    With \((s_1, s_2) = (1, 2)\), at \((z_1, z_2) = (1, 2)\), \(D = \Phi(1)\Phi(4) -
    \Phi(2)^2 = 0.84132 - 0.95502 = -0.11370 < 0\), while at \((z_1, z_2) = (-2,
    -1)\), \(D = \Phi(-2)^2 - \Phi(-1)\Phi(-4) = 0.000518 - 0.000005 = +0.000513 >
    0\). A kernel totally positive of order \(2\) has every \(2 \times 2\) minor
    nonnegative, so \(\Phi(sz)\) is not.
\end{proof}

The density kernel does not share the obstruction. Under \(u = z^2/2\),
\(\lambda = s^2\) the Gaussian density \(\phi(sz)\) becomes the Laplace kernel
\(e^{-\lambda u}\), sign-regular of every
order~\cite[Chapter~I]{karlin-totalpo-1968}. The variation-diminishing argument
of \cref{sec:S-upper} therefore works on this density kernel, through the
differentiated error~\eqref{eq:S-deriv}, and closes the upper bound there.

\section{The proof-dependency graph}
\label{sec:S-dag}

\Cref{tab:S-index} indexes the principal results. The setup
(\cref{lem:S-boundary}) supplies the boundary values that
\cref{thm:S-haar,lem:S-diff} consume. The Haar property (\cref{thm:S-haar}), via
the zero-counting lemma (\cref{lem:S-exp-sum-zeros}), feeds the lower bound
(\cref{cor:S-equiosc}) and, through the difference zero bound (\cref{lem:S-diff})
and the de la Vall\'ee Poussin proposition (\cref{prop:S-dlvp}), the uniqueness
chain (\cref{cor:S-global,cor:S-unique}). The tangent-Wronskian identity
(\cref{thm:S-wronskian}) is proved independently, by symbolic computation. The
upper bound (\cref{thm:S-upper}) is independent of the Haar route. It rests only
on the scale-mixture equation~\eqref{eq:S-stefanski} and the
double-peel zero count (\cref{rem:S-peel}) of the Gaussian density kernel. The
lower and upper bounds combine into the exact \(2K\) count on \((0, \infty)\),
which the symmetry lemma (\cref{lem:S-odd}) extends to \(4K\) on \(\R\). Existence
(\cref{prop:S-existence,prop:S-interior}) certifies that an interior minimiser
exists for \(K \le 11\). That branch is the graph's only apparent circularity.
\Cref{prop:S-existence} invokes \cref{prop:S-interior,lem:S-face} to exclude the
limiting faces, and \cref{prop:S-interior} itself rests on \cref{lem:S-face}. In
turn both assume the closure bound at level \(K - 1\), which is what
\cref{prop:S-existence} delivers one level down. The three are therefore resolved
together by induction on \(K\) from the \(K = 1\) base, as \cref{sec:S-existence}
sets out, and nothing in the chain is circular.

\begin{table}[htbp]
    \footnotesize
    \caption{Index of the appendices' principal results and their dependencies.}
    \label{tab:S-index}
    \centering
    \setlength{\tabcolsep}{4pt}
    \begin{tabular}{@{}l L{4.9cm} L{4.6cm}@{}}
        \toprule
        Result & Statement & Depends on \\
        \midrule
        \cref{lem:S-boundary}      & boundary values \(h(0) = h(\infty) = 0\) & setup \\
        \addlinespace[2pt]
        \cref{prop:S-existence}    & \raggedright attainment on \(\overline{\F_K}\), interior for
                                     \(1 \le K \le 11\) with error at least \(\Delta_K^\star\)
                                     throughout \(\overline{\F_K}\)
                                   & compactness, certified strict decrease, \cref{lem:S-face}, and
                                     \cref{prop:S-interior} \\
        \addlinespace[2pt]
        \cref{prop:S-interior}     & interior attainment (\(2 \le K \le 11\))
                                   & strict monotonicity, \cref{lem:S-face}, and
                                     \cref{prop:S-existence} at \(K - 1\) \\
        \addlinespace[2pt]
        \cref{lem:S-face}          & limiting faces exceed \(\Delta_K^\star\)
                                   & certified ratio bound and \cref{prop:S-existence} at \(K - 1\) \\
        \addlinespace[2pt]
        \cref{lem:S-exp-sum-zeros} & exponential-sum zero count \(\le 2K - 1\) & classical \\
        \addlinespace[2pt]
        \cref{thm:S-haar}          & \(V\) Haar of dimension \(2K - 1\)
                                   & \cref{lem:S-exp-sum-zeros,lem:S-boundary} \\
        \addlinespace[2pt]
        \cref{cor:S-equiosc}       & \(\ge 2K\) alternations & \cref{thm:S-haar} \\
        \addlinespace[2pt]
        \cref{thm:S-upper}         & \(\le 2K\) extrema, exact \(2K\)/\(4K\)
                                   & scale-mixture equation~\eqref{eq:S-stefanski},
                                     \cref{rem:S-peel},
                                     \cref{cor:S-equiosc}, and \cref{lem:S-odd} \\
        \addlinespace[2pt]
        \cref{rem:S-peel}          & \raggedright double-peel zero count \(\le 2K\) under
                                     an almost-everywhere positive density with finite
                                     polynomial moments
                                   & Rolle's theorem \\
        \addlinespace[2pt]
        \cref{lem:S-diff}          & difference has \(\le 2K - 2\) zeros
                                   & \cref{lem:S-exp-sum-zeros,lem:S-boundary} \\
        \addlinespace[2pt]
        \cref{prop:S-dlvp}         & de la Vall\'ee Poussin bound, equality case & \cref{lem:S-diff} \\
        \addlinespace[2pt]
        \cref{cor:S-global}        & global uniqueness & \cref{prop:S-dlvp,cor:S-equiosc} \\
        \addlinespace[2pt]
        \cref{cor:S-unique}        & strict local uniqueness & \cref{cor:S-global} \\
        \addlinespace[2pt]
        \cref{thm:S-wronskian}     & \(Q(z; s) \ne 0\) (\(K \le 5\)) & symbolic \\
        \addlinespace[2pt]
        \cref{conj:S-wronskian}    & identity for general \(K\) & verified \(K \le 12\) \\
        \addlinespace[2pt]
        \cref{lem:S-odd}           & error is odd, \(4K\) on \(\R\) & setup \\
        \addlinespace[2pt]
        \cref{lem:S-logconcave}    & \(\phi/\Phi\) strictly decreasing & elementary \\
        \addlinespace[2pt]
        \cref{prop:S-not-tp}       & \(\Phi(sz)\) not totally positive & numeric counterexample \\
        \bottomrule
    \end{tabular}
\end{table}

\section{Closed-form moments of the limiting mixing density}
\label{sec:S-moments}

The natural limiting object of the \(K\)-atom minimax family is the Kolmogorov--Smirnov mixing
density \(q\) of \cref{sec:intro}. Its low-order moments have closed forms, recorded here to
support \cref{sec:limit}. Write the scale-mixture
representation~\eqref{eq:S-stefanski} as \(Z = N S\), with \(Z \sim \mathrm{Logistic}(0,1)\),
\(N \sim N(0,1)\) independent, and \(\Lam(z) = \E[\Phi(z / S)]\), so that \(q\) is the law of the
reciprocal kernel scale \(1/S\). The logistic variable \(Z\) has the Mellin transform
\(\E[|Z|^s] = 2\,\Gamma(s+1)\,\eta_{\mathrm D}(s)\), where \(\eta_{\mathrm D}(s) = (1 - 2^{1-s})\,\zeta(s)\)
is the Dirichlet eta function. Dividing by the half-normal moments \(\E[|N|^s]\) and separating by parity
gives, for \(k \ge 1\),
\begin{equation}
    \begin{gathered}
        \E[S^{2k}] = 2^{k+1}\,k!\,\eta_{\mathrm D}(2k)
            = \frac{(4^k - 2)\,|B_{2k}|\,\pi^{2k}}{(2k-1)!!},\\
        \E[S^{2k+1}] = (2k+1)!!\,\eta_{\mathrm D}(2k+1)\,\sqrt{2\pi},
    \end{gathered}
    \label{eq:S-moments}
\end{equation}
so the even moments are rational multiples of \(\pi^{2k}\) through the Bernoulli numbers \(B_{2k}\) (for
instance \(\E[S^2] = \pi^2/3\)). The odd moments carry the odd zeta values (for instance
\(\E[S^3] = \tfrac94\sqrt{2\pi}\,\zeta(3)\)). Korolev and Shevtsova~\cite{korolev-moments-2026}
give the moments of the Kolmogorov distribution of all real orders in terms of the Riemann zeta and
Euler gamma functions, and with them the absolute moments of the logistic distribution of orders
\(p > -1\). Since \(S\) is twice a Kolmogorov variable in the
representation~\eqref{eq:S-stefanski}, the moments~\eqref{eq:S-moments} and their negative-order
counterparts are special cases of that result.

The kernel scale \(1/S\) inverts the moment formulas~\eqref{eq:S-moments} as
\(\E_q[s^k] = \E[S^{-k}]\), and its log-scale law has the closed-form mean and variance
\begin{equation}
\begin{aligned}
    \operatorname{mean}_q(\log s) &= \tfrac12\bigl(\gamma + \log 2 - 2\log\pi\bigr) = -0.50954\ldots,\\
    \operatorname{Var}_q(\log s) &= \frac{\pi^2}{24} - 2(\log 2)^2 - 4\log 2\,\log\pi - 2\,\zeta''(0)
        - \log^2(\pi/2) = 0.08524\ldots
\end{aligned}
    \label{eq:S-logvar}
\end{equation}
Here \(\gamma\) is the Euler--Mascheroni constant and \(\zeta''(0)\) is the
second derivative of the Riemann zeta function at the origin.

The inversion that produced the log-scale statistics~\eqref{eq:S-logvar} also reaches the
second moment of \(q\) through the continuation of the Mellin relation
\(\E[|Z|^s] = 2\,\Gamma(s+1)\,\eta_{\mathrm D}(s)\) to negative argument. At \(s = -2\) the zero
\(\eta_{\mathrm D}(-2) = 0\) cancels the pole of \(\Gamma(s + 1)\), and the surviving
derivative term, divided by the continued half-normal moment, evaluates through
\(\zeta'(-2) = -\zeta(3)/(4\pi^2)\) to
\(\eta_{\mathrm D}'(-2) = 7\zeta(3)/(4\pi^2)\), giving
\begin{equation}
    \E_q[s^2] = 2\,\eta_{\mathrm D}'(-2) = \frac{7\zeta(3)}{2\pi^2} = 0.42627\ldots,
    \label{eq:S-negmoment}
\end{equation}
the constant that the certified atoms' second moments approach (\cref{sec:limit}). The
shipped \(K = 10\) certificate reproduces the value~\eqref{eq:S-negmoment} to seven digits.

\section{A two-interval equilibrium reading of the plateau}
\label{sec:S-two-interval}

The finite plateau conjectured in \cref{sec:beyond} is the equality
\(\Delta_K^\star = \Delta_{11}^\star\) for \(K = 12, 13\), while at \(K = 14, \dots, 17\) mixtures
certified under the relaxed Bernstein--Akhiezer criterion carry audited supremum errors below
\(\Delta_{11}^\star\). The conjectured plateau has a potential-theoretic reading that reframes the open
plateau questions in terms of the limiting mixing density of \cref{sec:S-moments}, without resolving
them. Take the atom count \(K\) as an approximation
degree, the closure of the optimal
log-scale positions \(\{x_i^\star = \log s_i^\star\}\) as the support of the limit density
\(\rho_\infty\), and the asymptotic atom density as an equilibrium measure. On a single interval that
measure is the arcsine law and the minimax error falls monotonically with the degree. On two intervals
the extremal problem lives on a genus-one hyperelliptic curve, the equilibrium measure carries a
component on each interval, and the error can stall until the degree lets the alternation count split
across the two components~\cite{akhiezer-theorya-1956, widom-extremal-1969,
saff-logarith-2024}.

The plateau then reads as
single-interval saturation over \(K^\star \le K < K^{\star\star}\), with \(K^\star = 11\) the last
degree carrying a strict certificate. The certified \(K = 14, 15\) improvements move the onset of
falling to \(K^{\star\star} = 14\), and the support edges do not single out one
bifurcation degree. The upper edge of the best-known log-scale support lies at \(+0.48\) for the
\(K = 14\) witness and \(+0.69\) for the \(K = 16\) certificate. The \(K = 15\)
witness carries an atom of weight \(5.6 \times 10^{-10}\) at \(+0.78\), with its
remaining mass ending at \(+0.48\).
Whether the two-interval mechanism, mass migrating past the single-interval edge as the degree grows,
is what ends the conjectured \(K = 12, 13\) plateau remains open.

This reading of the open questions is not a proof. It carries a terminological caution, since
``Bernstein--Akhiezer'' names both the relaxed certificate criterion of \cref{sec:beyond} and the classical
two-interval equilibrium measure, and only the latter is meant here.

\section{Method comparison in causal mediation}
\label{sec:S-comparison}

This appendix gives the full comparison summarised in \cref{sec:qb-comparison},
with the measured costs of the pointwise route inside the sampling loop that
\cref{sec:pointwise} lists.

\Cref{tab:comparison} compares three inner-integral strategies across three
outcome-prevalence regimes at \(R = 10^4\) outer draws. Its caption records the
agreement among the certified mixture, inner Monte Carlo, and the pointwise
evaluator, and the wall times, and \cref{sec:qb-comparison} draws the
conclusions. The repeated-sampling experiment behind the coverage figures there
is the script \texttt{inst/reproduction/Coverage.R} in the repository.

The same draws give the odds ratio \eqref{eq:exact-or} of the natural indirect
effect. Across the low-, moderate-, and high-prevalence regimes of
\cref{tab:comparison}, the mixture point estimates and \(95\%\) intervals are
\(6.27\) \((4.50, 8.46)\), \(7.59\) \((4.89, 10.77)\), and \(13.99\)
\((6.48, 22.10)\), against population values of \(5.27\), \(5.72\), and
\(6.41\). The closed form \eqref{eq:rare-form} gives \(18.8\) \((8.23, 44.6)\),
\(24.1\) \((7.71, 78.3)\), and \(225\) \((11.3, 4438)\) on the same data sets,
which have \(v = 2.56\) at the population coefficients.

The pointwise rows carry one qualification. At high prevalence the evaluator
returns a nonfinite value on three of the \(10^4\) draws at truncation order
\(10\) and on four at order \(30\), which \cref{tab:comparison} marks, and those
draws are dropped before the quantiles are taken. The account below traces
the failures to the two edges of the evaluator's domain.

\begin{table}[htbp]
    \footnotesize
    \caption{Method comparison for the natural indirect effect (NIE) and its
        95\% quasi-Bayesian (QB) confidence interval (CI) at \(R = 10^4\) outer
        draws (\(n = 500\)) across three outcome-prevalence regimes. The minimax
        mixture, inner Monte Carlo (MC), and Pirjol's pointwise Poisson-summation
        evaluator at truncation orders \(10\) and \(30\) evaluate the inner integral
        accurately and agree throughout, with the mixture and the pointwise route agreeing
        to all six decimals the script carries in the low- and moderate-prevalence
        regimes. The high-prevalence aggregates differ in the fourth decimal, since the
        pointwise quantiles there are taken over a reduced draw set. A dagger marks
        a pointwise entry with nonfinite draws dropped, three of \(10^4\) at order \(10\)
        and four at order \(30\), all in the high-prevalence regime. Six fall below the
        evaluator's lower domain edge, where the sampled total mediator coefficient comes
        closest to zero, and the remaining order-\(30\) draw crosses that order's upper
        edge. Reproducible via \texttt{inst/reproduction/Example1.R}, whose sweep also
        runs the repository default \(K = 22\), which agrees with the displayed \(K = 11\)
        values to the shown precision. The pointwise rows need the upstream sources,
        which are third-party and not shipped, and the script runs without them. The
        parameter draws are formed from a Cholesky factor, which is unique for a
        positive definite covariance, so the sweep reproduces across BLAS and LAPACK
        builds. That reproducibility was verified by running the sweep under OpenBLAS and
        under a reference BLAS, where all \(45\) numeric records agree exactly. Every wall
        time reported in this paper was measured on a single 16-core workstation.}
    \label{tab:comparison}
    \centering
    \begin{tabular}{lllrrl}
        \toprule
        Scenario & Method & \(\widehat{\mathrm{NIE}}_R\) & 95\% QB CI & Bias & Wall time \\
        \midrule
        low (\(P_Y = 0.38\)) & Minimax \(K = 11\) & 0.4112
                             & \((0.3491, 0.4695)\) & \(-0.004\) & \(\le 0.01\) s \\
                             & Inner MC, \(M = 10^4\) & 0.4112
                             & \((0.3488, 0.4695)\) & \(-0.004\) & \(4.8\) s \\
                             & Pirjol, \(N = 10\) & 0.4112
                             & \((0.3491, 0.4695)\) & \(-0.004\) & \(0.17\) s \\
                             & Pirjol, \(N = 30\) & 0.4112
                             & \((0.3491, 0.4695)\) & \(-0.004\) & \(0.43\) s \\
        \addlinespace
        moderate (\(P_Y = 0.69\)) & Minimax \(K = 11\) & 0.3178
                             & \((0.2227, 0.4004)\) & \(-0.005\) & \(\le 0.01\) s \\
                             & Inner MC, \(M = 10^4\) & 0.3178
                             & \((0.2230, 0.4006)\) & \(-0.005\) & \(4.8\) s \\
                             & Pirjol, \(N = 10\) & 0.3178
                             & \((0.2227, 0.4004)\) & \(-0.005\) & \(0.17\) s \\
                             & Pirjol, \(N = 30\) & 0.3178
                             & \((0.2227, 0.4004)\) & \(-0.005\) & \(0.43\) s \\
        \addlinespace
        high (\(P_Y = 0.82\)) & Minimax \(K = 11\) & 0.2438
                             & \((0.1436, 0.3311)\) & \(-0.004\) & \(\le 0.01\) s \\
                             & Inner MC, \(M = 10^4\) & 0.2438
                             & \((0.1438, 0.3308)\) & \(-0.004\) & \(4.8\) s \\
                             & Pirjol, \(N = 10^{\dagger}\) & 0.2439
                             & \((0.1438, 0.3311)\) & \(-0.004\) & \(0.17\) s \\
                             & Pirjol, \(N = 30^{\dagger}\) & 0.2439
                             & \((0.1437, 0.3311)\) & \(-0.004\) & \(0.43\) s \\
        \bottomrule
    \end{tabular}
\end{table}

Within the mediation workflow the pointwise route is a usable inner
engine, and \cref{tab:comparison} measures it rather than assuming it, and each
route corroborates the other where the two agree. Three costs remain. First,
the truncation order must be managed per draw, so the loop-wide error statement is
empirical rather than a certified constant fixed in advance. Second, derivatives
for delta-method intervals and sensitivity analyses require numerical
differentiation of a procedural routine. Third,
and specific to the sampling loop, the evaluator returns a nonfinite value at
both ends of its domain. The lower edge lies near \(v = 0.0265\) and does not move
with the truncation order. The upper edge depends on the mean of the inner normal
as well as on \(v\) at both orders. At order \(10\) it lies near \(v = 142\) at a
zero mean, above every draw of this regime, while at order \(30\) it lies within
the sampled range, so no single threshold in \(v\) describes the failures. The
pointwise route is undefined on three of the \(10^4\) high-prevalence draws of
\cref{tab:comparison} at order \(10\) and on four at order \(30\). All three
failures at order \(10\), and three of the four at order \(30\), fall below the
lower edge, where the total mediator coefficient \(\theta_2 + \theta_3 z\) is
near zero. Its smallest sampled magnitude in this regime is \(0.075\), against
\(1.10\) and \(0.76\) in the other two. The remaining order-\(30\) failure is
the largest sampled variance of the regime, \(v = 56.5\), against maxima of \(9.8\)
and \(13.2\) elsewhere. Raising the order for accuracy therefore
narrows the domain from above, while the order-independent lower edge lies
where the indirect effect is near null, and an unguarded loop propagates the first
such draw into the interval.

\end{appendices}

\bibliographystyle{unsrtnat}
\bibliography{references-numalg}

\end{document}